\documentstyle[12pt]{article}

\font\twlgot =eufm10 scaled \magstep1 \font\egtgot =eufm8
\font\sevgot =eufm7 \font\twlmsb =msbm10 scaled \magstep1
\font\egtmsb =msbm8 \font\sevmsb =msbm7

\newfam\gotfam
\def\pgot{\fam\gotfam\twlgot}
\textfont\gotfam\twlgot \scriptfont\gotfam\egtgot
\scriptscriptfont\gotfam\sevgot
\def\got{\protect\pgot}
\newfam\msbfam
\textfont\msbfam\twlmsb \scriptfont\msbfam\egtmsb
\scriptscriptfont\msbfam\sevmsb
\def\Bbb{\protect\pBbb}
\def\pBbb{\relax\ifmmode\expandafter\Bb\else\typeout{You cann't use
Bbb in text mode}\fi}
\def\Bb #1{{\fam\msbfam\relax#1}}

\newcommand{\gQ}{{\got T}}

\newcommand{\gA}{{\got A}}

\newcommand{\gd}{{\got d}}

\newcommand{\gS}{{\got S}}

\def\thebibliography#1{\bigskip\bigskip\centerline{\bf References}\list
  {[\arabic{enumi}]}{\settowidth\labelwidth{#1}\leftmargin\labelwidth
    \advance\leftmargin\labelsep
    \usecounter{enumi}}
    \def\newblock{\hskip .11em plus .33em minus .07em}
    \sloppy\clubpenalty4000\widowpenalty4000
    \sfcode`\.=1000\relax}

\def\op#1{\mathop{\fam0 #1}\limits}

\newcommand{\im}{{\rm Im\,}}
\newcommand{\nm}[1]{|{#1}|}
\newcommand{\beq}{\begin{equation}}
\newcommand{\eeq}{\end{equation}}
\newcommand{\ben}{\begin{eqnarray}}
\newcommand{\een}{\end{eqnarray}}
\newcommand{\be}{\begin{eqnarray*}}
\newcommand{\ee}{\end{eqnarray*}}
\newcommand{\bea}{\begin{eqalph}}
\newcommand{\eea}{\end{eqalph}}
\newcommand{\cA}{{\cal A}}

\newcommand{\cP}{{\cal P}}

\newcommand{\cL}{{\cal L}}

\newcommand{\cE}{{\cal E}}

\newcommand{\cS}{{\cal S}}
\newcommand{\cC}{{\cal C}}
\newcommand{\cO}{{\cal O}}

\newcommand{\bL}{{\bf L}}

\newcommand{\bE}{{\bf E}}

\newcommand{\al}{\alpha}
\newcommand{\vr}{\varrho}
\newcommand{\bt}{\beta}
\newcommand{\dl}{\delta}
\newcommand{\la}{\lambda}
\newcommand{\La}{\Lambda}
\newcommand{\f}{\phi}
\newcommand{\om}{\omega}

\newcommand{\m}{\mu}

\newcommand{\G}{\Gamma}
\newcommand{\th}{\theta}

\newcommand{\vt}{\vartheta}

\newcommand{\up}{\upsilon}

\newcommand{\di}{{\rm dim\,}}

\newcommand{\si}{\sigma}
\newcommand{\Si}{\Sigma}
\newcommand{\w}{\wedge}

\newcommand{\ol}{\overline}
\newcommand{\dr}{\partial}
\newcommand{\ar}{\op\longrightarrow}
\newcommand{\ot}{\otimes}

\newcommand{\e}{\epsilon}

\newcommand{\rdr}{\stackrel{\leftarrow}{\dr}{}}
\newcommand{\lto}{\leftarrow}
\newcommand{\llr}{\op\longleftarrow}

\newcounter{remark}
\newcounter{example}
\newcounter{theorem}
\newcounter{proposition}
\newcounter{lemma}
\newcounter{corollary}
\newcounter{definition}
\def\theremark{\arabic{remark}}
\def\thetheorem{\arabic{theorem}}
\def\theproposition{\arabic{proposition}}
\def\thelemma{\arabic{lemma}}
\def\thecorollary{\arabic{corollary}}
\def\thedefinition{\arabic{definition}}

\def\theremark{\arabic{remark}}
\def\thetheorem{\arabic{theorem}}
\def\theproposition{\arabic{proposition}}
\def\thelemma{\arabic{lemma}}
\def\thecorollary{\arabic{corollary}}
\def\thedefinition{\arabic{definition}}

\newenvironment{proof}{{\it Proof.}}{\hfill q.e.d.\medskip}
\newenvironment{rem}{\refstepcounter{remark}\medskip{\bf
Remark \theremark.}}{\medskip}

\newenvironment{theo}{\refstepcounter{theorem}\medskip
{\bf Theorem \thetheorem.} \it}{\medskip}
\newenvironment{prop}{\refstepcounter{proposition}\medskip
{\bf Proposition \theproposition.}\it}{\medskip}
\newenvironment{lem}{\refstepcounter{lemma}\medskip
{\bf Lemma \thelemma.}\it}{\medskip}
\newenvironment{cor}{\refstepcounter{corollary}\medskip
{\bf Corollary \thecorollary.}\it}{\medskip}
\newenvironment{defi}{\refstepcounter{definition}\medskip
{\bf Definition \thedefinition.}\it}{\medskip}

\newcommand{\mar}[1]{}

\begin{document}
\hbox{}



\centerline{\bf NOETHER'S INVERSE SECOND THEOREM IN HOMOLOGY
TERMS}

\bigskip

\centerline{\sc Giovanni Giachetta, Luigi Mangiarotti, Gennadi
Sardanashvily}
\bigskip
\bigskip

{\small \centerline{\bf Abstract}
\bigskip

 A generic degenerate Lagrangian system of even and odd
variables on an arbitrary smooth manifold is examined in terms of
the Grassmann-graded variational bicomplex. Its Euler--Lagrange
operator obeys Noether identities which need not be independent,
but satisfy first-stage Noether identities, and so on. However,
non-trivial higher-stage Noether identities are ill defined,
unless a certain homology condition holds. We show that, under
this condition, there exists the exact Koszul--Tate chain complex
whose boundary operator produces all non-trivial Noether and
higher-stage Noether identities of an original Lagrangian system.
Noether's inverse second theorem that we prove associates to this
complex a cochain sequence whose ascent operator provides all
gauge and higher-stage gauge supersymmetries of an original
Lagrangian. }

\bigskip
\bigskip
\bigskip
\centerline{\bf Introduction}
\bigskip

Since Noether identities of a Lagrangian systems of even variables
are parameterized by elements of a Grassmann algebra, we address
from the beginning a generic degenerate Lagrangian system of even
and odd variables on an arbitrary smooth manifold. It is described
in terms of the Grassmann-graded variational bicomplex
\cite{barn,jmp05,cmp04} generalizing the well-known variational
bicomplex for even Lagrangian systems on fiber bundles
\cite{ander,jmp,tak2} (Section 1). Theorem \ref{v11} provides its
relevant cohomology.

Any Euler--Lagrange operator obeys trivial Noether identities
which are defined as boundaries of a certain chain complex
(Definition \ref{v120}). A Lagrangian system is said to be
degenerate if its Euler--Lagrange operator obeys non-trivial
Noether identities given by homology of this complex. Noether
identities need not be independent, but satisfy non-trivial
first-stage Noether identities, which in turn are subject to the
second-stage ones, and so on. Thus, we have a hierarchy of
reducible Noether identities. A problem is that trivial
higher-stage Noether identities need not be boundaries.

The notion of reducible Noether identities has come from that of
reducible constraints. Their Koszul--Tate complex has been
invented by analogy with that of constraints under a rather
restrictive regularity condition that field equations as well as
Noether identities of arbitrary stage can be locally separated
into the independent and dependent ones \cite{fisch}. This
condition has also come from the case of constraints locally given
by a finite number of functions which the inverse mapping theorem
is applied to. A problem is that, in contrast with constraints,
Noether and higher-stage Noether identities are differential
operators. They are locally given by a set of functions and their
jet prolongations on an infinite order jet manifold. Since the
latter is a Fr\'echet, but not Banach manifold, the inverse
mapping theorem fails to be valid.

In Section 2, we show that, if Noether and higher-stage Noether
identities are finitely generated and iff a certain homology
regularity condition (Definition \ref{v155}) holds, one can
associate to the Euler--Lagrange operator of a degenerate
Grassmann-graded Lagrangian system the exact Koszul--Tate complex
(\ref{w3}) whose boundary operator (\ref{w4}) produces all
non-trivial Noether and higher-stage Noether identities  (Theorem
\ref{v163}).

Noether's second theorems in different formulations relate the
Noether and higher-stage Noether identities to the gauge and
higher-stage gauge symmetries and supersymmetries of a Lagrangian
system \cite{jpa05,jmp05,fulp}. In Section 3, we prove Noether's
inverse second theorem (Theorem \ref{w35}) which associates to the
above mentioned Koszul--Tate complex the cochain sequence
(\ref{w36}), whose ascent operator (\ref{w108}) provides gauge and
higher-stage gauge supersymmetries of an original Lagrangian
system. This operator need not be nilpotent. Therefore, a
formulation of Noether's direct second theorem in cohomology terms
meets difficulties. However, the ascent operator admits a
nilpotent extension and the above mentioned cochain sequence is a
complex if gauge and higher-stage gauge supersymmetries of an
original Lagrangian system constitute an algebra (Remark
\ref{gg6}).

In Section 4, an example of a reducible degenerate Lagrangian
system coming from the topological BF theory is examined in
detail.

\bigskip
\bigskip
\centerline{\bf 1. Preliminary. Grassmann-graded Lagrangian
systems}
\bigskip

Smooth manifolds throughout are real, finite-dimensional,
Hausdorff, second-countable (hence, paracompact) and connected.
The symbols $\La$, $\Si$, $\Xi$ stand for symmetric multi-indices,
e.g., $\La=(\la_1...\la_k)$, $|\La|=k$, and
$\la+\La=(\la\la_1...\la_k)$.

Let $Y\to X$, $\di X=n$, be a fiber bundle and $J^rY$, $r\in \Bbb
N$, the jet manifolds of its sections. The index $r=0$ stands for
$Y$. There is the inverse system of affine bundles
\mar{5.10}\beq
X\op\longleftarrow^\pi Y\op\longleftarrow^{\pi^1_0} J^1Y
\longleftarrow \cdots J^{r-1}Y \op\longleftarrow^{\pi^r_{r-1}}
J^rY\longleftarrow\cdots, \label{5.10}
\eeq
whose projective limit $(J^\infty Y;\pi^\infty_r:J^\infty Y\to
J^rY)$ is a paracompact Fr\'echet manifold \cite{tak2}. A bundle
atlas $(U_Y;x^\la,y^i)$ of $Y\to X$ induces the coordinate atlas
\mar{jet1}\beq
((\pi^\infty_0)^{-1}(U_Y); x^\la, y^i, y^i_\la,\ldots,
y^i_\La,\ldots), \qquad 0\leq|\La|, \label{jet1}
\eeq
of $J^\infty Y$. The inverse system (\ref{5.10}) yields the direct
system
\mar{5.7}\beq
\cO^*X\op\longrightarrow^{\pi^*} \cO^*Y
\op\longrightarrow^{\pi^1_0{}^*} \cO_1^*Y \ar\cdots \cO^*_{r-1}Y
\op\longrightarrow^{\pi^r_{r-1}{}^*}
 \cO_r^*Y \longrightarrow\cdots  \label{5.7}
\eeq
of graded differential algebras (henceforth GDAs) $\cO_r^*Y$ of
exterior forms on jet manifolds $J^rY$ with respect to the
pull-back monomorphisms $\pi^r_{r-1}{}^*$. Its direct limit is the
GDA  $\cO_\infty^*Y$ of all exterior forms on finite order jet
manifolds modulo the pull-back identification. The GDA
$\cO_\infty^*Y$ is split into the above mentioned variational
bicomplex describing Lagrangian systems on a fiber bundle $Y\to
X$.

\begin{rem} \label{gg30} \mar{gg30}
The GDA $\cO_\infty^*Y$ is a subalgebra of the GDA considered in
\cite{tak2}. Let $\gS^*_r$ be the sheaf of germs of exterior forms
on $J^rY$ and $\ol\gS^*_r$ its canonical presheaf (we follow the
terminology of \cite{hir}). There is the direct system of
presheaves
\be
 \ol\gS^*\ar \ol\gS^*_1 \ar\cdots
\ol\gS^*_r \ar\cdots,
\ee
whose direct limit is a presheaf of GDAs on $J^\infty Y$. Let
$\gQ^*_\infty $ be the sheaf of germs of this presheaf. The
structure module $\G\gQ^*_\infty$ of sections of this sheaf is a
GDA such that, given an element $\f\in \G\gQ^*_\infty$, there
exist an open neighbourhood $U$ of each point of $J^\infty Y$ and
an exterior form $\f^{(k)}$ on some finite order jet manifold
$J^kY$ so that $\f|_U= \pi^{\infty*}_k\f^{(k)}|_U$. There is an
obvious monomorphism $\cO^*_\infty Y\to \G\gQ^*_\infty$. Note that
$J^*_\infty Y$ admits the partition of unity by elements of the
ring $\G\gQ^0_\infty$, but not $\cO_\infty^0Y$. Therefore, one can
obtain cohomology of $\G\gQ^*_\infty Y$ by virtue the abstract de
Rham theorem \cite{ander,tak2}. The GDA $\cO_\infty^*Y$ is proved
to possesses the same cohomology as $\G\gQ^*_\infty$
\cite{lmp,jmp}. The above mentioned Theorem \ref{v11} is similarly
proved.
\end{rem}

In order to describe Noether identities generated by elements of
projective Grassmann-graded $C^\infty(X)$-modules of finite rank,
we appeal to the well-known Serre--Swan theorem, extended to
noncompact manifolds \cite{book05,ren}, and to its following
combination \cite{jmp05a,book05} with the Batchelor theorem
\cite{bart}.

\begin{prop} \label{v0} \mar{v0}
Given a smooth manifold $Z$, the  exterior algebra of a projective
$C^\infty(Z)$-module of finite rank is isomorphic to the ring of
graded functions on some graded manifold whose body is $Z$.
\end{prop}

Let $(Z,\gA)$ be a graded manifold with a body $Z$ and a structure
sheaf $\gA$ of Grassmann $C^\infty_Z$-algebras of finite rank,
where $C^\infty_Z$ is the sheaf of germs of smooth real functions
on $Z$ \cite{bart}. The above mentioned Batchelor theorem states
an isomorphism of $(Z,\gA)$ to a graded manifold $(Z,\gA_Q)$ with
the structure sheaf $\gA_Q$ of germs of sections of an exterior
bundle
\be
\w Q^*=\Bbb R\op\oplus_Z Q^*\op\oplus_Z\op\w^2
Q^*\op\oplus_Z\cdots,
\ee
where $Q^*$ is the dual of some vector bundle $Q\to Z$. In our
case, Batchelor's isomorphism is fixed from the beginning. Let us
call $(Z,\gA_Q)$ a graded manifold modelled over $Q$. Its
structure ring $\cA_Q$ of graded functions consists of global
sections of the exterior bundle $\w Q^*$. Let $\gd\cA_Q$ be the
real Lie superalgebra $\gd\cA_Q$ of graded derivations of the
$\Bbb R$-ring $\cA_Q$, i.e.,
\be
u(ff')=u(f)f'+(-1)^{[u][f]}fu(f'), \qquad f,f'\in \cA_Q, \qquad
u\in \gd\cA_Q,
\ee
where the symbol $[.]$ stands for the Grassmann parity. Then the
Chevalley--Eilenberg complex of $\gd\cA_Q$ with coefficients in
$\cA_Q$ can be constructed \cite{fuks}.  Its subcomplex
$\cS^*[Q;Z]$ of $\cA_Q$-linear morphisms is the Grassmann-graded
Chevalley--Eilenberg differential calculus
\be
0\to \Bbb R\to \cA_Q \ar^d \cS^1[Q;Z]\ar^d\cdots
\cS^k[Q;Z]\ar^d\cdots
\ee
over $\cA_Q=\cS^0[Q;Z]$. The graded exterior product $\w$ and the
Chevalley--Eilenberg coboundary operator $d$ make $\cS^*[Q;Z]$
into a bigraded differential algebra (henceforth BGDA)
\mar{v22}\beq
\f\w\f' =(-1)^{|\f||\f'| +[\f][\f']}\f'\w \f, \qquad  d(\f\w\f')=
d\f\w\f' +(-1)^{|\f|}\f\w d\f',  \label{v22}
\eeq
where $|.|$ denotes the form degree. Moreover, $\cS^*[Q;Z]$ is a
minimal differential calculus over $\cA_Q$ generated by elements
$df$, $f\in \cA_Q$. There is a natural monomorphism $\cO^*Z\to
\cS^*[Q;Z]$.

Elements of the BGDA $\cS^*[Q;Z]$ can be seen as graded exterior
forms on a manifold $Z$. Given an open subset $U\subset Z$, let
$\cA_U$ be the Grassmann algebra of sections of the sheaf $\gA_Q$
over $U$, and let $\cS^*[Q;U]$ be the Chevalley--Eilenberg
differential calculus over $\cA_U$. Given an open subset
$U'\subset U$, the restriction morphism $\cA_U\to\cA_{U'}$ yields
a homomorphism of BGDAs $\cS^*[Q;U]\to \cS^*[Q;U']$. Thus,  we
obtain the presheaf $\{U,\cS^*[Q;U]\}$ of BGDAs on a manifold $Z$
and the sheaf $\gQ^*[Q;Z]$ of germs of this presheaf. Moreover,
$\{U,\cS^*[Q;U]\}$ is a canonical presheaf of $\gQ^*[Q;Z]$. Hence,
$\cS^*[Q;Z]$ is the BGDA of global sections of the sheaf
$\gQ^*[Q;Z]$, and there is the restriction morphism $\cS^*[Q;Z]\to
\cS^*[Q;U]$ for any open $U\subset Z$. Due to this morphism,
elements of $\cS^*[Q;Z]$ take the following local form. Bundle
coordinates $(z^A,q^a)$ on $Q$ and the corresponding fiber basis
$\{c^a\}$ for $Q^*$ provide a local basis $(z^A, c^a)$ for the
graded manifold $(Z,\gA_Q)$ such that graded functions on $Z$ read
\mar{v23}\beq
f=\op\sum_{0\leq k} \frac1{k!}f_{a_1\ldots a_k}c^{a_1}\cdots
c^{a_k}, \qquad f_{a_1\ldots a_k}\in C^\infty(Z). \label{v23}
\eeq
Owing to the isomorphism $VQ= Q\times Q$, the fiber basis
$\{\dr_a\}$ for the vertical tangent bundle $VQ\to Q$ of $Q\to Z$
is the dual of $\{c^a\}$. Then the $\cA_Q$-module $\gd\cA_Q$ of
graded derivations is locally generated by the elements $\dr_A$,
$\dr_a$ acting on graded functions (\ref{v23}) by the rule
\be
\dr_A(f)=\op\sum_{0\leq k} \frac1{k!}\dr_A(f_{a_1\ldots
a_k})c^{a_1}\cdots c^{a_k}, \quad \dr_d(f)= \op\sum_{0\leq k}
\frac1{k!}\op\sum_{1\leq i\leq k} (-1)^{i-1}f_{a_1\ldots
a_k}c^{a_1}\cdots \dl^{a_i}_d\cdots c^{a_k}.
\ee
Relative to the dual bases $\{dz^A\}$ for $T^*Z$ and $\{dc^b\}$
for $Q^*$, the BGDA $\cS^*[Q;Z]$ is locally generated by graded
one-forms $dz^A$, $dc^a$.

A generic Lagrangian system of even and odd variables on a smooth
manifold $X$ is defined in terms of composite graded manifolds
whose bodies are a fiber bundle $Y\to X$ and its jet manifolds
$J^rY$ \cite{jmp05a} (see \cite{jmp05,cmp04} for a particular case
of an affine bundle $Y\to X$). Let $F\to X$ be a vector bundle.
Let us consider the graded manifold $(J^rY,, \gA_{F_r})$ modelled
over the product $F_r=J^rY\times_XJ^rF$. There is an epimorphism
of graded manifolds $(J^{r+1}Y,\gA_{F_{r+1}}) \to
(J^rY,\gA_{F_r})$. It consists of the surjection $\pi^{r+1}_r$ and
the sheaf monomorphism $\pi_r^{r+1*}\gA_{F_r}\to \gA_{F_{r+1}}$,
where $\pi_r^{r+1*}\gA_{F_r}$ is the pull-back of the topological
fiber bundle $\gA_{F_r}\to J^rY$ onto $J^{r+1}Y$. This
monomorphism of sheaves yields a monomorphism of their canonical
presheaves $\ol \gA_{F_r}\to \ol \gA_{F_{r+1}}$ which associates
to every open subset $U\subset J^{r+1}Y$ the ring of sections of
$\gA_{F_r}$ over $\pi^{r+1}_r(U)$. Accordingly, there is a
monomorphism of graded commutative rings $\cA_{F_r} \to
\cA_{F_{r+1}}$ which induces the monomorphism of BGDAs
\mar{v4}\beq
\cS^*[F_r;J^rY]\to \cS^*[F_{r+1};J^{r+1}Y]. \label{v4}
\eeq
As a consequence, we have the direct system of BGDAs
\mar{v6}\beq
\cS^*[Y\op\times_X F;Y]\ar \cS^*[F_1;J^1Y]\ar\cdots
\cS^*[F_r;J^rY]\ar\cdots. \label{v6}
\eeq
Its direct limit $\cS^*_\infty[F;Y]$  is a BGDA of all graded
exterior forms on jet manifolds $J^rY$ modulo monomorphisms
(\ref{v4}). The relations (\ref{v22}) hold. Monomorphisms
$\cO^*_rY\to \cS^*[F_r;J^rY]$ provide a monomorphism of the direct
system (\ref{5.7}) to the direct system (\ref{v6}) and, thus, a
monomorphism $\cO^*_\infty Y\to \cS^*_\infty[F;Y]$ of their direct
limits. Moreover, $\cS^*_\infty[F;Y]$ is an $\cO^0_\infty
Y$-algebra.

Elements of the BGDA $\cS^*_\infty[F;Y]$ can be seen as graded
exterior forms on $J^\infty Y$. Indeed, let $\gS^*[F_r;J^rY]$ be
the sheaf of BGDAs on $J^rY$ and $\ol\gS^*[F_r;J^rY]$ its
canonical presheaf whose elements are the Chevalley--Eilenberg
differential calculus over elements of the presheaf
$\ol\gA_{F_r}$. Then the presheaf monomorphisms $\ol \gA_{F_r}\to
\ol \gA_{F_{r+1}}$ yield the direct system of presheaves
\mar{v15}\beq
\ol\gS^*[Y\times F;Y]\ar \ol\gS^*[F_1;J^1Y] \ar\cdots
\ol\gS^*[F_r;J^rY]  \ar\cdots. \label{v15}
\eeq
Its direct limit is a presheaf of BGDAs on $J^\infty Y$. Let
$\gQ^*_\infty[F;Y]$ be the sheaf of germs of this presheaf. The
structure module $\G\gQ^*_\infty[F;Y]$ of sections of
$\gQ^*_\infty[F;Y]$ is a BGDA such that, given an element $\f\in
\G\gQ^*_\infty[F;Y]$, there exist an open neighbourhood $U$ of
each point of $J^\infty Y$ and a graded exterior form $\f^{(k)}$
on some finite order jet manifold $J^kY$ so that $\f|_U=
\pi^{\infty*}_k\f^{(k)}|_U$. There is a monomorphism
$\cS^*_\infty[F;Y] \to\G\gQ^*_\infty[F;Y]$ (cf. that in Remark
\ref{gg30}).

Due to this monomorphism, one can restrict $\cS^*_\infty[F;Y]$ to
the coordinate chart (\ref{jet1}), and say that
$\cS^*_\infty[F;Y]$ as an $\cO^0_\infty Y$-algebra is locally
generated by  the elements
\be
(c^a_\La, dx^\la,\th^a_\La=dc^a_\La-c^a_{\la+\La}dx^\la,\th^i_\La=
dy^i_\La-y^i_{\la+\La}dx^\la), \qquad 0\leq |\La|.
\ee
One calls $(y^i,c^a)$ the local basis for $\cS^*_\infty[F;Y]$. We
further use the collective symbol $s^A$ for its elements, together
with the notation $s^A_\La$, $\th^A_\La=ds^A_\La-
s^A_{\la+\La}dx^\la$, and $[A]=[s^A]$. Let $\gd \cS^0_\infty[F;Y]$
be the Lie superalgebra of graded derivations of the $\Bbb R$-ring
$\cS^0_\infty[F;Y]$. Its elements read
\mar{gg}\beq
 \vt=\vt^\la\dr_\la + \op\sum_{0\leq|\La|} \vt_\La^A\dr^\La_A,
 \qquad \vt^\la, \vt_\La^A
 \in \cS^0[F;Y], \label{gg}
\eeq
where $\dr^\La_A(s_\Si^B)=\dl_A^B\dl^\La_\Si$ up to permutations
of multi-indices $\La$ and $\Si$.
 The interior product
$\vt\rfloor\f$ and the Lie derivative $\bL_\vt\f$,
$\f\in\cS^*_\infty[F;Y]$, obey the relations
\be
&& \vt\rfloor \f=\vt^\la\f_\la +
\op\sum_{0\leq|\La|}(-1)^{[\f_A^\La][A]}\vt^A_\La \f_A^\La, \qquad
\f\in \cS^1_\infty[F;Y],\\
&& \vt\rfloor(\f\w\si)=(\vt\rfloor \f)\w\si
+(-1)^{|\f|+[\f][\vt]}\f\w(\vt\rfloor\si), \qquad \f,\si\in
\cS^*_\infty[F;Y], \\
&& \bL_\vt\f=\vt\rfloor d\f+ d(\vt\rfloor\f), \qquad
\bL_\vt(\f\w\si)=\bL_\vt(\f)\w\si
+(-1)^{[\vt][\f]}\f\w\bL_\vt(\si).
\ee
In particular, the total derivatives are defined as graded
derivations
\be
\gd \cS^0_\infty[F;Y]\ni d_\la =\dr_\la + \op\sum_{0\leq|\La|}
s_{\la+\La}^A\dr^\La_A, \qquad  d_\la\f=\bL_{d_\la}\f, \qquad
d_\La=d_{\la_1}\cdots d_{\la_k}.
\ee

The BGDA $\cS^*_\infty[F;Y]$ is split into
$\cS^0_\infty[F;Y]$-modules $\cS^{k,r}_\infty[F;Y]=h_k\circ
h^r(\cS^*_\infty[F;Y])$ of $k$-contact and $r$-horizontal graded
forms. Accordingly, the graded exterior differential $d$ on
$\cS^*_\infty[F;Y]$ falls into the sum $d=d_H+d_V$ of the total
and vertical differentials where
\be
d_H\circ h_k=h_k\circ d\circ h_k,
 \qquad d_H(\f)= dx^\la\w d_\la(\f).
\ee
These differentials together with the graded projection
endomorphism
\be
\vr=\op\sum_{k>0} \frac1k\ol\vr\circ h_k\circ h^n, \qquad
\ol\vr(\f)= \op\sum_{0\leq|\La|} (-1)^{\nm\La}\th^A\w
[d_\La(\dr^\La_A\rfloor\f)], \qquad \f\in \cS^{>0,n}_\infty[F;Y],
\ee
and the graded variational operator $\dl=\vr\circ d$ make the BGDA
$\cS^*_\infty[F;Y]$ into the above mentioned Grassmann-graded
variational bicomplex. We restrict our consideration to its short
variational subcomplex and the subcomplex of one-contact graded
forms
\mar{g111,2}\ben
&& 0\to \Bbb R\ar \cS^0_\infty[F;Y]\ar^{d_H}\cS^{0,1}_\infty[F;Y]
\cdots \ar^{d_H} \cS^{0,n}_\infty[F;Y]\ar^\dl \bE_1
=\vr(\cS^{1,n}_\infty[F;Y]), \label{g111}\\
&& 0\to \cS^{1,0}_\infty[F;Y]\ar^{d_H} \cS^{1,1}_\infty[F;Y]\cdots
\ar^{d_H}\cS^{1,n}_\infty[F;Y]\ar^\vr \bE_1\to 0. \label{g112}
\een
One can think of their even elements
\mar{0709,'}\ben
&& L=\cL\om\in \cS^{0,n}_\infty[F;Y], \qquad \om=dx^1\w\cdots \w
dx^n,
\label{0709}\\
&& \dl L= \th^A\w \cE_A\om=\op\sum_{0\leq|\La|}
 (-1)^{|\La|}\th^A\w d_\La (\dr^\La_A L)\om\in \bE_1 \label{0709'}
\een
as being a graded Lagrangian and its Euler--Lagrange operator,
respectively.

\begin{theo} \label{v11} \mar{v11}
(i) Cohomology of the complex (\ref{g111}) equals the de Rham
cohomology $H^*(Y)$ of $Y$. (ii) The complex (\ref{g112}) is exact.
\end{theo}

\begin{proof} The proof follows that of \cite{cmp04},
Theorem 2.1 (see Section 5).
\end{proof}

\begin{cor} \label{cmp26} \mar{cmp26}
A $\dl$-closed (i.e., variationally trivial) graded density
$L\in \cS^{0,n}_\infty[F;Y]$ reads
\mar{g215}\beq
L=h_0\psi + d_H\xi, \qquad \xi\in \cS^{0,n-1}_\infty[F;Y],
\label{g215}
\eeq
where $\psi$ is a closed $n$-form on $Y$. In particular, a
$\dl$-closed odd graded density is $d_H$-exact.
\end{cor}

\begin{cor} \label{cmp26'} \mar{cmp26'}
Any graded density $L$ admits the decomposition
\mar{g99}\beq
dL=\dl L - d_H\Xi,
\qquad \Xi\in \cS^{1,n-1}_\infty[F;Y], \label{g99}\\
\eeq
where $L+\Xi$ is a Lepagean equivalent of $L$.
\end{cor}

A graded derivation $\vt$ (\ref{gg}) is called contact if the Lie
derivative $\bL_\vt$ preserves the ideal of contact graded forms
of the BGDA $\cS^*_\infty[F;Y]$. It reads
\mar{gg1}\beq
\vt=\vt_H +\vt_V= [\vt^\la d_\la] +[\up^A\dr_A + \op\sum_{0<|\La|}
d_\La\up^A\dr_A^\La], \qquad \up^A=\vt^A- s^A_\la \vt^\la.
\label{gg1}
\eeq
A contact graded derivation $\vt$ (\ref{gg1}) is called a
variational supersymmetry of a Lagrangian $L$ (\ref{0709}) if the
Lie derivative $\bL_\vt L$ is $d_H$-exact. The following holds
\cite{cmp04}.

\begin{prop}
A contact graded derivation $\vt$ (\ref{gg1}) is a variational
supersymmetry of $L$ iff its vertical part $\vt_V$, vanishing on
$\cO^*X\subset \cS^*_\infty[F;Y]$, is well.
\end{prop}

Therefore, we further restrict our consideration to vertical
contact graded derivations
\mar{0672}\beq
\vt=\up^A\dr_A + \op\sum_{0<|\La|} d_\La\up^A\dr_A^\La.
\label{0672}
\eeq
Such a derivation is completely defined by its first summand
$\up=\up^A\dr_A$.

\begin{prop}
As a result of the splitting (\ref{g99}), the Lie derivative
$\bL_\vt L$ of a Lagrangian $L$ along a vertical contact graded
derivations $\vt$ (\ref{0672}) admits the decomposition
\mar{g107}\beq
\bL_\vt L= \up\rfloor\dl L +d_H(\vt\rfloor \Xi)). \label{g107}
\eeq
\end{prop}

\begin{prop}
An odd vertical contact graded derivations $\vt$ (\ref{0672}) is a
variational supersymmetry of $L$ iff the odd density
$\up\rfloor\dl L=\up^A\cE_A\om$ is $d_H$-exact.
\end{prop}

A vertical contact graded derivation $\vt$ (\ref{0672}) is called
nilpotent if $\bL_\vt(\bL_\vt\f)=0$ for any horizontal graded form
$\f\in \cS^{0,*}_\infty[F;Y]$. It is nilpotent only if it is odd
and iff
\mar{0688}\beq
\vt(\up)=\vt(\up^A\dr_A)=\op\sum_{0\leq|\Si|}
\up^B_\Si\dr^\Si_B(\up^A)\dr_A=0. \label{0688}
\eeq

For the sake of simplicity, the common symbol $\up$ further stands
for $\vt$ (\ref{0672}), its summand $\up$, and the Lie derivative
$\bL_\vt$. We agree to call $\up$ the graded derivation of the
BGDA $\cS^*_\infty[F;Y]$.

\begin{rem} \label{w74} \mar{w74}
Right contact graded derivations $\op\up^\lto
={\op\dr^\lto}_A\up^A$ of the BGDA $\cS^*_\infty[F;Y]$ are also
involved in the sequel. They act on graded forms $\f$ on the right
by the rule
\be
\op\up^\lto(\f)=\op d^\lto(\f)\lfloor \op\up^\lto +\op
d^\lto(\f\lfloor\op\up^\lto), \qquad
\op\up^\lto(\f\w\f')=(-1)^{[\f'][\op\up^\lto]}\op\up^\lto(\f)\w\f'+
\f\w\op\up^\lto(\f').
\ee
For instance, ${\op\dr^\lto}_A(\f)=(-1)^{([\f]+1)[A]}\dr_A(\f)$,
${\op d^\lto}_\La=d_\La$ and ${\op d^\lto}_H(\f)=
(-1)^{|\f|}d_H(\f)$. With right graded derivations, we have the
right Euler--Lagrange operator
\be
\op\dl^\lto L= {\op\cE^\lto}_A\om\w \th^A, \qquad {\op\cE^\lto}_A
=\op\sum_{0\leq|\La|}
 (-1)^{|\La|}d_\La (\rdr^\La_A (L)).
\ee
An odd right graded derivation $\op\up^\lto$
is a variational supersymmetry of a graded Lagrangian $L$ iff the odd
graded  density
${\op\cE^\lto}_A\up^A\om$ is $d_H$-exact.
\end{rem}

\begin{rem} Any
local graded functions $f'$, $f^\La$, $0\leq |\La|\leq k$, and a
graded exterior form $\f$ obey the equalities
\mar{0606a-c}\ben
&& \op\sum_{0\leq |\La|\leq k} f^\La d_\La f'\om= \op\sum_{0\leq
|\La|\leq k} (-1)^{|\La|}d_\La (f^\La) f'\om + d_H\si,
\label{0606a}\\
&& \op\sum_{0\leq |\La|\leq k} (-1)^{|\La|}d_\La(f^\La \f)=
\op\sum_{0\leq |\La|\leq k} \eta (f)^\La d_\La \f, \qquad
(\eta\circ\eta)(f)^\La=f^\La, \label{0606b}
\\ && \eta (f)^\La = \op\sum_{0\leq|\Si|\leq
k-|\La|}(-1)^{|\Si+\La|} C^{|\Si|}_{|\Si+\La|} d_\Si f^{\Si+\La},
\qquad C^a_b=\frac{b!}{a!(b-a)!}. \label{0606c}
\een
In particular, the decomposition (\ref{g107}) takes the local form
(\ref{0606a}), but Corollary \ref{cmp26'} states that the second
term in its right-hand side is globally $d_H$-exact.
\end{rem}

\bigskip
\bigskip
\centerline{\bf 2. The Koszul--Tate complex of Noether identities}
\bigskip

Given a degenerate Grassmann-graded Lagrangian system
$(\cS^*_\infty[F;Y],L)$, let us associate to the Euler--Lagrange
operator $\dl L$ (\ref{0709'}) of a graded Lagrangian $L$
(\ref{0709}) the exact Koszul--Tate chain complex with the
boundary operator whose nilpotency conditions provide all
non-trivial Noether and higher-stage Noether identities for $\dl
L$.

\begin{rem} We introduce the following notation. Let $E\to X$ be a vector bundle
and $E^*$ its dual. The bundle product
\be
\ol E^*=E^*\op\ot_X\op\w^n T^*X
\ee
is called the density-dual of $E$. Given the pull-back $E_Y$ of
$E$ onto $Y$, let us consider the BGDA $\cS^*_\infty[F;E_Y]$.
There are monomorphisms of $\cO^0_\infty Y$-algebras
$\cS^*_\infty[F;Y]\to \cS^*_\infty[F;E_Y]$ and $\cO^*_\infty E\to
\cS^*_\infty[F;E_Y]$ whose images contain the common subalgebra
$\cO^*_\infty Y$. We consider: (i) the subring $\cP^0_\infty
E_Y\subset \cO^0_\infty E_Y$ of polynomial functions in fiber
coordinates  of the vector bundles $J^rE_Y\to J^rY$, (ii) the
corresponding subring
$\cP^0_\infty[F;E_Y]\subset\cS^0_\infty[F;E_Y]$ of graded
functions with polynomial coefficients belonging to $\cP^0_\infty
E_Y$, (iii) the  subalgebra  $\cP^*_\infty[F;Y;E]$ of the BGDA
$\cS^*_\infty[F;E_Y]$ over the subring $\cP^0_\infty[F;E_Y]$.
Given vector bundles $V,V',E,E'$ over $X$, let us denote
\mar{v90}\beq
\cP^*_\infty[V'V;F;Y;EE']= \cP^*_\infty[V'\op\times_X
V\op\times_XF;Y;E\op\times_X E']. \label{v90}
\eeq
The BGDA $\cP^*_\infty[F;Y;E]$ and, similarly, the BGDA
(\ref{v90}) possess the same  cohomology as $\cS^*_\infty[F;Y]$ in
Theorem \ref{v11}. Since $H^*(Y)=H^*(E_Y)$, this cohomology of the
BGDA $\cS^*_\infty[F;Y]$ equals that of the BGDA
$\cS^*_\infty[F;E_Y]$. Furthermore, one can replace the BGDA
$\cS^*_\infty[F;E_Y]$ with $\cP^*_\infty[F;Y;E]$ in the condition
of Theorem \ref{v11} due to the fact that sheaves of $\cP^0_\infty
E_Y$-modules are also sheaves of $\cO^0_\infty Y$-modules.
\end{rem}

\begin{rem}
For the sake of simplicity, we assume that the vertical tangent
bundle $VY$ of a fiber bundle $Y\to X$ admits the splitting
$VY=Y\times_X W$, where $W\to X$ is some vector bundle. In this
case, there no fiber bundles under consideration whose transition
functions can vanish on the shell $\cE_A=0$. Let $\ol Y^*$ denote
the density-dual of $W$ in this splitting.
\end{rem}

Let us enlarge the BGDA $\cS^*_\infty[F;Y]$ to the BGDA
$\cP^*_\infty[\ol Y^*;F;Y;\ol F^*]$ whose local basis is
\be
\{s^A, \ol s_A\}, \qquad [\ol s_A]=([A]+1){\rm mod}\,2.
\ee
Following the physical terminology \cite{barn}, we agree to call
$\ol s_A$ the antifields of antifield number Ant$[\ol s_A]= 1$.
The BGDA $\cP^*_\infty[\ol Y^*;F;Y;\ol F^*]$ is provided with the
nilpotent right graded derivation
\mar{gg33}\beq
\ol\dl=\rdr^A \cE_A, \label{gg33}
\eeq
where $\cE_A$ are the graded variational derivatives
(\ref{0709'}). We call $\ol\dl$ the Koszul--Tate differential.

\begin{defi} One  says  that an
element of the BGDA $\cP^*_\infty[\ol Y^*;F;Y;\ol F^*]$ or its
extension vanishes on the shell if it is $\ol\dl$-exact.
\end{defi}

With the Koszul--Tate differential (\ref{gg33}), the module
$\cP^{0,n}_\infty[\ol Y^*;F;Y;\ol F^*]$ of graded densities is
split into the chain complex
\be
 0\lto \cS^{0,n}_\infty[F;Y]
\llr^{\ol\dl} \cP^{0,n}_\infty[\ol Y^*;F;Y;\ol F^*]_1\cdots
\llr^{\ol\dl} \cP^{0,n}_\infty[\ol Y^*;F;Y;\ol F^*]_k \cdots
\ee
graded by the antifield number. Let us consider its subcomplex
\mar{v042}\beq
0\lto \im\ol\dl \llr^{\ol\dl} \cP^{0,n}_\infty[\ol Y^*;F;Y;\ol
F^*]_1 \llr^{\ol\dl} \cP^{0,n}_\infty[\ol Y^*;F;Y;\ol F^*]_2.
\label{v042}
\eeq
It is exact at $\im\ol\dl$. Let us examine its first homology
$H_1(\ol\dl)$.

\begin{rem}
If there is no danger of confusion, elements of homology of a
chain complex are identified to its representatives. A chain
complex is called $r$-exact if its homology of degree $k\leq r$ is
trivial.
\end{rem}

A generic one-chain of the complex (\ref{v042}) takes the form
\mar{0712}\beq
\Phi= \op\sum_{0\leq|\La|} \Phi^{A,\La}\ol s_{\La A} \om, \qquad
\Phi^{A,\La}\in \cS^0_\infty[F;Y]. \label{0712}
\eeq
The cycle condition reads
\mar{0713}\beq
\ol\dl \Phi=\op\sum_{0\leq|\La|} \Phi^{A,\La} d_\La \cE_A \om=0.
\label{0713}
\eeq
This equality is a Noether identity which the graded variational
derivatives $\cE_A$ (\ref{0709'}) satisfy. Conversely, any
equality of the form (\ref{0713}) comes from some cycle
(\ref{0712}). A Noether identity (\ref{0713}) is trivial  if a
cycle is a boundary
\be
 \Phi=
\op\sum_{0\leq|\La|,|\Si|} T^{(A\La)(B\Si)}d_\Si\cE_B\ol s_{\La
A}\om, \qquad T^{(A\La)(B\Si)}=-(-1)^{[A][B]} T^{(B\Si)(A\La)}.
\ee

\begin{defi} \label{v120} \mar{v120}
Noether identities which the Euler--Lagrange operator $\dl L$
(\ref{0709'}) satisfies are one-cycles of the chain complex
(\ref{v042}). Trivial Noether identities are boundaries.
Non-trivial Noether identities considered modulo the trivial ones
are non-zero elements of the first homology $H_1(\ol\dl)$ of the
chain complex (\ref{v042}).
\end{defi}

One can say something more if the $\cS^0_\infty[F;Y]$-module
$H_1(\ol \dl)$ is finitely generated. Namely, there exists a
projective Grassmann-graded $C^\infty(X)$-module $\cC_{(0)}\subset
H_1(\ol \dl)$ of finite rank such that any element $\Phi\in
H_1(\ol \dl)$ factorizes via elements of $\cC_{(0)}$ as
\mar{v63,71}\ben
&& \Phi= \op\sum_{0\leq|\Xi|} G^{r,\Xi} d_\Xi \Delta_r\om, \qquad
G^{r,\Xi}\in
\cS^0_\infty[F;Y], \label{v63}\\
&&\Delta_r=\op\sum_{0\leq|\La|} \Delta_r^{A,\La}\ol s_{\La A},
\qquad \Delta_r^{A,\La}\in \cS^0_\infty[F;Y], \label{v71}
\een
where $\{\Delta_r\}$ is a local basis for $\cC_{(0)}$. This means
that any Noether identity (\ref{0713}) is a corollary of the
Noether identities
\mar{v64}\beq
\ol\dl\Delta_r= \op\sum_{0\leq|\La|} \Delta_r^{A,\La} d_\La
\cE_A=0. \label{v64}
\eeq
Clearly, the factorization (\ref{v63}) is independent of
specification of local bases $\{\Delta_r\}$. By virtue of the
Serre--Swan theorem, the module $\cC_{(0)}$ is isomorphic to a
module of sections of the product $\ol V^*\times_X \ol E^*$, where
$\ol V^*$ and $\ol E^*$ are the density-duals of some vector
bundles $V\to X$ and $E\to X$.

\begin{defi}
If the first homology $H_1(\ol \dl)$ of the chain complex
(\ref{v042}) is finitely generated, its generating elements
$\Delta_r\in \cC_{(0)}$ (\ref{v71}) and the corresponding
equalities (\ref{v64}) are called the complete Noether identities.
\end{defi}

For instance, let $L$ (\ref{0709}) be a variationally trivial
Lagrangian. Its Euler--Lagrange operator $\dl L=0$ obeys the
Noether identities which are finitely generated by the Noether
identities $\Delta_A=\ol s_A$.

\begin{prop} \label{v137} \mar{v137}
If the homology $H_1(\ol\dl)$ of the chain complex (\ref{v042}) is
finitely generated, this complex  can be extended to the one-exact
chain complex (\ref{v66}) with a boundary operator whose
nilpotency conditions are equivalent to the complete Noether
identities (\ref{v64}).
\end{prop}

\begin{proof}
Let us enlarge the BGDA $\cP^*_\infty[\ol Y^*;F;Y;\ol F^*]$ to  the
BGDA
\mar{w1}\beq
\cP^*_\infty[\ol E^*\ol Y^*;F;Y;\ol F^*\ol V^*], \label{w1}
\eeq
possessing the local basis
\be
\{s^A,\ol s_A, \ol c_r\}, \qquad [\ol c_r]=([\Delta_r]+1){\rm
mod}\,2, \qquad {\rm Ant}[\ol c_r]=2.
\ee
The BGDA (\ref{w1}) is provided with the nilpotent right graded
derivation
\be
\dl_0=\ol\dl + \rdr^r\Delta_r,
\ee
called the zero-stage Koszul--Tate differential. It is readily
observed that its nilpotency conditions (\ref{0688}) are
equivalent to the complete Noether identities (\ref{v64}). Then
the module $\cP^{0,n}_\infty[\ol E^*\ol Y^*;F;Y;\ol F^*\ol
V^*]_{\leq 3}$ of graded densities of antifield number
Ant$[\f]\leq 3$ is split into the chain complex
\mar{v66}\ben
&&0\lto \im\ol\dl \llr^{\ol\dl} \cP^{0,n}_\infty[\ol Y^*;F;Y;\ol
F^*]_1\llr^{\dl_0}
\cP^{0,n}_\infty[\ol E^*\ol Y^*;F;Y;\ol F^*\ol V^*]_2 \label{v66}\\
&& \qquad \llr^{\dl_0} \cP^{0,n}_\infty[\ol E^*\ol Y^*;F;Y;\ol
F^*\ol V^*]_3. \nonumber
\een
Let $H_*(\dl_0)$ denote its homology. We have
$H_0(\dl_0)=H_0(\ol\dl)=0$. Furthermore, any one-cycle $\Phi$ up
to a boundary takes the form (\ref{v63}) and, therefore, it is a
$\dl_0$-boundary
\be
\Phi= \op\sum_{0\leq|\Si|} G^{r,\Xi} d_\Xi \Delta_r\om
=\dl_0(\op\sum_{0\leq|\Si|} G^{r,\Xi}\ol c_{\Xi r}\om).
\ee
Hence, $H_1(\dl_0)=0$, i.e., the complex (\ref{v66}) is one-exact.
\end{proof}

Let us examine the second homology $H_2(\dl_0)$ of the complex
(\ref{v66}). A generic two-chain  reads
\mar{v77}\beq
\Phi= G + H= \op\sum_{0\leq|\La|} G^{r,\La}\ol c_{\La r}\om +
\op\sum_{0\leq|\La|,|\Si|} H^{(A,\La)(B,\Si)}\ol s_{\La A}\ol
s_{\Si B}\om.
\label{v77}
\eeq
The cycle condition takes the form
\mar{v79}\beq
 \dl_0 \Phi=\op\sum_{0\leq|\La|} G^{r,\La}d_\La\Delta_r\om +\ol\dl H=0.
\label{v79}
\eeq
This is a first-stage Noether identity which the complete Noether
identities (\ref{v71}) satisfy. Conversely, let
\be
\Phi=\op\sum_{0\leq|\La|} G^{r,\La}\ol c_{\La r}\om\in
\cP^{0,n}_\infty[\ol E^*\ol Y^*;F;Y;\ol F^*\ol V^*]_2
\ee
be a graded density such that the first-stage Noether identity
(\ref{v79}) hold. This identity is obviously a cycle condition of
the two-chain (\ref{v77}).

\begin{defi} \label{gg2} \mar{gg2} The first-stage Noether identities which the
complete Noether identities satisfy are two-cocycles of the
one-exact chain complex (\ref{v66}).
\end{defi}

The first-stage Noether identity (\ref{v79}) is trivial either if
a two-cycle $\Phi$ (\ref{v77}) is a boundary or its summand $G$,
linear in antifields, vanishes on the shell. Because of the second
requirement, trivia first-stage Noether identities need not be
two-boundaries, unless the following condition is satisfied.

\begin{defi} \label{gg3} \mar{gg3} One says that the chain
complex (\ref{v66}) obeys the
two-homology regularity condition if any $\ol\dl$-cycle $\f\in
\cP^{0,n}_\infty[\ol Y^*;F;Y;\ol F^*]_2$ is a $\dl_0$-boundary.
\end{defi}

\begin{prop} \label{v134} \mar{v134}
Non-trivial first-stage Noether identities are identified to
non-zero elements of the second homology $H_2(\dl_0)$ of the
complex (\ref{v66}) iff the two-homology regularity condition
hold.
\end{prop}

\begin{proof}
It suffices to show that, if the summand $G$ of a two-cycle $\Phi$
(\ref{v77}) is $\ol\dl$-exact, $\Phi$ is a boundary. If $G=\ol\dl
\Psi$, then
\mar{v169}\beq
\Phi=\dl_0\Psi +(\ol \dl-\dl_0)\Psi + H. \label{v169}
\eeq
The cycle condition reads
\be
\dl_0\Phi=\ol\dl((\ol\dl-\dl_0)\Psi + H)=0.
\ee
Then $(\ol \dl-\dl_0)\Psi + H$ is $\dl_0$-exact since any
$\ol\dl$-cycle $\f\in \cP^{0,n}_\infty[\ol Y^*;F;Y;\ol F^*]_2$, by
assumption, is a $\dl_0$-boundary. Consequently, $\Phi$
(\ref{v169}) is $\dl_0$-exact. Conversely, let $\Phi\in
\cP^{0,n}_\infty[\ol Y^*;F;Y;\ol F^*]_2$ be an arbitrary
$\ol\dl$-cycle.
 The cycle condition reads
\mar{v100}\beq
\ol\dl\Phi= 2\Phi^{(A,\La)(B,\Sigma)}\ol s_{\La A} \ol\dl\ol
s_{\Sigma B}\om= 2\Phi^{(A,\La)(B,\Sigma)}\ol s_{\La A} d_\Si
\cE_B\om=0. \label{v100}
\eeq
It follows that $\Phi^{(A,\La)(B,\Sigma)} \ol\dl\ol s_{\Sigma
B}=0$ for all indices $(A,\La)$. We obtain
\be
\Phi^{(A,\La)(B,\Sigma)} \ol s_{\Sigma B}= G^{(A,\La)(r,\Xi)}d_\Xi
\Delta_r +\ol\dl S^{(A,\La)}.
\ee
Hence, $\Phi$ takes the form
\mar{v135}\beq
\Phi=G'^{(A,\La)(r,\Xi)} d_\Xi\Delta_r \ol s_{\La A}\om +\ol\dl
S^{(A,\La)}\ol s_{\La A}\om. \label{v135}
\eeq
We can associate to $\Phi$ (\ref{v135}) the three-chain
\be
\Psi= G'^{(A,\La)(r,\Xi)} \ol c_{\Xi r} \ol s_{\La A}\om +
S^{(A,\La)}\ol s_{\La A}\om
\ee
such that
\be
\dl_0\Psi=\Phi +\si = \Phi + G''^{(A,\La)(r,\Xi)}d_\La\cE_A \ol
c_{\Xi r} \om + S'^{(A,\La)}\ol\dl\ol s_{\La A}\om.
\ee
Owing to the equality $\ol\dl\Phi=0$, we have $\dl_0\si=0$. Since
the term $G''$ of $\si$ is $\ol\dl$-exact, then $\si$  by
assumption is $\dl_0$-exact, i.e., $\si=\dl_0\psi$. It follow that
$\Phi=\dl_0\Psi -\dl_0\psi$.
\end{proof}

If the two-homology regularity condition is satisfied, let us
suppose that the second homology $H_2(\dl_0)$ of the complex
(\ref{v66}) is finitely generated as follows. There exists a
projective Grassmann-graded $C^\infty(X)$-module $\cC_{(1)}\subset
H_2(\dl_0)$ of finite rank such that any element $\Phi\in
H_2(\dl_0)$ factorizes via elements of $\cC_{(1)}$ as
\mar{v80,1}\ben
&& \Phi= \op\sum_{0\leq|\Xi|} \Phi^{r_1,\Xi} d_\Xi
\Delta_{r_1}\om, \qquad \Phi^{r_1,\Xi}\in
\cS^0_\infty[F;Y], \label{v80}\\
&&\Delta_{r_1}=G_{r_1}+ h_{r_1}=\op\sum_{0\leq|\La|}
\Delta_{r_1}^{r,\La}\ol c_{\La r} + h_{r_1}, \qquad
 h_{r_1}\om\in
\cP^{0,n}_\infty[\ol Y^*;F;Y;\ol F^*], \label{v81}
\een
where $\{\Delta_{r_1}\}$ is a local basis for $\cC_{(1)}$. Thus,
any first-stage Noether identity (\ref{v79}) results from the
equalities
\mar{v82}\beq
 \op\sum_{0\leq|\La|} \Delta_{r_1}^{r,\La} d_\La \Delta_r +\ol\dl
h_{r_1} =0. \label{v82}
\eeq
By virtue of the Serre--Swan theorem, $\cC_{(1)}$ is isomorphic to
the module of sections of the product $\ol V^*_1\op\times_X \ol
E^*_1$, where $\ol V^*_1$ and $\ol E^*_1$ are the density-duals of
some vector bundles $V_1\to X$ and $E_1\to X$.

\begin{defi} \mar{gg4} \label{gg4}
(i) If the chain complex (\ref{v66}) obeys the two-homology
regularity condition and its second homology $H_2(\dl_0)$ is
finitely generated, the generating elements $\Delta_{r_1}\in
\cC_{(1)}$ (\ref{v81}) of $H_2(\dl_0)$ and the corresponding
equalities (\ref{v82}) are called the complete first-stage Noether
identities. (ii) A degenerate Lagrangian system is said to be
one-stage reducible if it possesses complete Noether and
first-stage Noether identities.
\end{defi}

In other words, a degenerate Lagrangian system
$(\cS^*_\infty[F;Y],L)$ is first-stage reducible if one associates
to it a one-exact chain complex (\ref{v66}) which obeys the
two-homology regularity condition and whose second homology is
finitely generated.

\begin{prop} \label{v139} \mar{v139} The one-exact chain complex (\ref{v66})
associated to a first-stage reducible degenerate Lagrangian system
can be extended to the two-exact chain complex (\ref{v87}) with a
boundary operator whose nilpotency conditions are equivalent to
complete Noether and first-stage Noether identities.
\end{prop}

\begin{proof}
Let us consider the BGDA $\cP^*_\infty[\ol E^*_1\ol E^*\ol
Y^*;F;Y;\ol F^*\ol V^*\ol V^*_1]$ possessing the local basis
\be
\{s^A,\ol s_A, \ol c_r, \ol c_{r_1}\}, \qquad [\ol
c_{r_1}]=([\Delta_{r_1}]+1){\rm mod}\,2 \qquad  {\rm Ant}[\ol
c_{r_1}]=3.
\ee
It can be provided with the nilpotent graded derivation
\be
\dl_1=\dl_0 + \rdr^{r_1} \Delta_{r_1},
\ee
called the first-stage Koszul--Tate differential. Ii is easily
seen that its nilpotency conditions (\ref{0688}) are equivalent to
the complete Noether identities (\ref{v64}) and complete
first-stage Noether identities (\ref{v82}). Then the module
$\cP^{0,n}_\infty[\ol E^*_1\ol E^*\ol Y^*;F;Y;\ol F^*\ol V^*\ol
V^*_1]_{\leq 4}$ of graded densities of antifield number
Ant$[\f]\leq 4$ is split into the chain complex
\mar{v87}\ben
&&0\lto \im\ol\dl \llr^{\ol\dl} \cP^{0,n}_\infty[\ol Y^*;F;Y;\ol
F^*]_1\llr^{\dl_0} \cP^{0,n}_\infty[\ol E^*\ol Y^*;F;Y;\ol F^*\ol
V^*]_2\llr^{\dl_1}
\label{v87}\\
&& \qquad \cP^{0,n}_\infty[\ol E^*_1\ol E^*\ol Y^*;F;Y;\ol F^*\ol
V^*\ol V^*_1]_3
 \llr^{\dl_1}
\cP^{0,n}_\infty[\ol E^*_1\ol E^*\ol Y^*;F;Y;\ol F^*\ol V^*\ol
V^*_1]_4. \nonumber
\een
Let $H_*(\dl_1)$ denote its homology. It is readily observed that
\be
H_0(\dl_1)=H_0(\ol\dl), \qquad H_1(\dl_1)=H_1(\dl_0)=0.
\ee
By virtue of the expression (\ref{v80}), any two-cycle of the
complex (\ref{v87}) is a boundary
\be
 \Phi= \op\sum_{0\leq|\Xi|} \Phi^{r_1,\Xi} d_\Xi \Delta_{r_1}\om
=\dl_1(\op\sum_{0\leq|\Xi|} \Phi^{r_1,\Xi} \ol c_{\Xi r_1})\om.
\ee
It follows that $H_2(\dl_1)=0$, i.e., the complex (\ref{v87}) is
two-exact.
\end{proof}

If the third homology $H_3(\dl_1)$ of the chain complex
(\ref{v87}) is not trivial, there are non-trivial second-stage
Noether identities which the first-stage ones satisfy, and so on.
Iterating the arguments, we come to the following.

Given a first-stage reducible degenerate Lagrangian system
$(\cS^*_\infty[F;Y],L)$ in accordance with Definition \ref{gg4},
let us assume the following.

(a) Given an integer $N\geq 1$, there are vector bundles
$V_1,\ldots, V_N, E_1, \ldots, E_N$ over $X$, and the BGDA
$\cS^*_\infty[F;Y]$ is enlarged to the BGDA
\mar{v91}\beq
\ol\cP^*_\infty\{N\}=\cP^*_\infty[\ol E^*_N\cdots\ol E^*_1\ol
E^*\ol Y^*;F;Y;\ol F^*\ol V^*\ol V^*_1\cdots\ol V_N^*] \label{v91}
\eeq
with a local basis $\{s^A,\ol s_A, \ol c_r, \ol c_{r_1}, \ldots,
\ol c_{r_N}\}$ graded by antifield numbers Ant$[\ol c_{r_k}]=k+2$.
Let the indexes $k=-1,0$ further stand for $\ol s_A$ and $\ol
c_r$, respectively.

(b) The BGDA $\ol\cP^*_\infty\{N\}$ (\ref{v91}) is provided with a
nilpotent graded derivation
\mar{v92,'}\ben
&&\dl_N=\rdr^A\cE_A +
\op\sum_{0\leq|\La|}\rdr^r\Delta_r^{A,\La}\ol s_{\La A} +
\op\sum_{1\leq k\leq N}\rdr^{r_k} \Delta_{r_k},
\label{v92}\\
&& \Delta_{r_k}=G_{r_k} + h_{r_k}= \op\sum_{0\leq|\La|}
\Delta_{r_k}^{r_{k-1},\La}\ol c_{\La r_{k-1}} + \op\sum_{0\leq
|\Si|, |\Xi|}(h_{r_k}^{(r_{k-2},\Si)(A,\Xi)}\ol c_{\Si
r_{k-2}}\ol s_{\Xi A}+...), \label{v92'}
\een
of antifield number -1.

(c) With $\dl_N$, the module $\ol\cP^{0,n}_\infty\{N\}_{\leq N+3}$
of graded densities of antifield number Ant$[\f]\leq N+3$ is split
into the $(N+1)$-exact chain complex
\mar{v94}\ben
&&0\lto \im \ol\dl \llr^{\ol\dl} \cP^{0,n}_\infty[\ol Y^*;F;Y;\ol
F^*]_1\llr^{\dl_0} \ol\cP^{0,n}_\infty\{0\}_2\llr^{\dl_1}
\ol\cP^{0,n}_\infty\{1\}_3\cdots
\label{v94}\\
&& \qquad
 \llr^{\dl_{N-1}} \ol\cP^{0,n}_\infty\{N-1\}_{N+1}
\llr^{\dl_N} \ol\cP^{0,n}_\infty\{N\}_{N+2}\llr^{\dl_N}
\ol\cP^{0,n}_\infty\{N\}_{N+3}, \nonumber
\een
which satisfies the following $(N+1)$-homology regularity
condition.

\begin{defi} \label{v155} \mar{v155} One says that the chain complex (\ref{v94})
obeys the $(N+1)$-homology regularity condition if any
$\dl_{k<N-1}$-cycle $\f\in \ol\cP_\infty^{0,n}\{k\}_{k+3}\subset
\ol\cP_\infty^{0,n}\{k+1\}_{k+3}$ is a $\dl_{k+1}$-boundary.
\end{defi}

Note that
the $(N+1)$-exactness of the complex (\ref{v94}) implies that any
$\dl_{k<N-1}$-cycle $\f\in \ol\cP_\infty^{0,n}\{k\}_{k+3}$, $k<N$,
is a $\dl_{k+2}$-boundary, but not necessary a
$\dl_{k+1}$-boundary.

If $N=1$, the complex $\ol\cP^{0,n}_\infty\{1\}_{\leq 4}$
(\ref{v94}) is the chain complex (\ref{v87}). Therefore, we agree
to call $\dl_N$ (\ref{v92}) the $N$-stage Koszul--Tate
differential. Its nilpotency implies the complete Noether
identities (\ref{v64}), first-stage Noether identities
(\ref{v82}), and the complete $(k\leq N)$-stage Noether identities
\mar{v93}\beq
\op\sum_{0\leq|\La|} \Delta_{r_k}^{r_{k-1},\La}d_\La
(\op\sum_{0\leq|\Si|} \Delta_{r_{k-1}}^{r_{k-2},\Si}\ol c_{\Si
r_{k-2}}) + \ol\dl(\op\sum_{0\leq |\Si|,
|\Xi|}h_{r_k}^{(r_{k-2},\Si)(A,\Xi)}\ol c_{\Si r_{k-2}}\ol s_{\Xi
A})=0,  \label{v93}
\eeq
which the complete $(k-1)$-stage Noether identities
$\Delta_{r_{k-1}}$ (\ref{v92'}) satisfy.

\begin{defi} \mar{gg5} \label{gg5} If the above mentioned
assumptions (a) -- (c) hold, a degenerate Grassmann-graded
Lagrangian system $(\cS^*_\infty[F;Y],L)$ is called $N$-stage
reducible.
\end{defi}

If the $(N+2)$-homology of the complex (\ref{v94}) is not trivial,
an $N$-stage reducible Lagrangian system is $(N+1)$-stage
reducible under the following conditions.

\begin{theo} \label{v163} \mar{v163} Given an $N$-stage reducible
Lagrangian system in accordance with Definition \ref{gg5}, let us
suppose that the $(N+2)$-homology $H_{N+2}(\dl_N)$ of the
associated chain complex (\ref{v94}) is not trivial. Then the
following holds.

(i)  The $(N+1)$-stage Noether identities which the complete
$N$-stage Noether identities satisfy are the $(N+2)$-cycles of the
complex (\ref{v94}), and vice versa.

(ii)  The trivial $(N+1)$-stage Noether identities are
$(N+2)$-boundaries iff the $(N+2)$-homology regularity condition
holds. In this case, non-trivial $(N+1)$-stage Noether identities
modulo the trivial ones are identified to non-zero elements of the
homology $H_{N+2}(\dl_N)$.

(iii) If the homology $H_{N+2}(\dl_N)$ is finitely generated, the
complex (\ref{v94}) admits an $(N+2)$-exact extension. The
nilpotency of its boundary operator implies the complete Noether
and $(k\leq N+1)$-stage Noether identities.
\end{theo}

\begin{proof} (i) A
generic $(N+2)$-chain $\Phi\in \ol\cP^{0,n}_\infty\{N\}_{N+2}$
takes the form
\mar{v156}\beq
\Phi= G + H= \op\sum_{0\leq|\La|} G^{r_N,\La}\ol c_{\La r_N}\om +
\op\sum_{0\leq |\Si|, |\Xi|}(H^{(A,\Xi)(r_{N-1},\Si)}\ol s_{\Xi
A}\ol c_{\Si r_{N-1}}+...)\om. \label{v156}
\eeq
The cycle condition $\dl_N\Phi=0$ implies the equality
\mar{v145}\beq
\op\sum_{0\leq|\La|} G^{r_N,\La}d_\La (\op\sum_{0\leq|\Si|}
\Delta_{r_N}^{r_{N-1},\Si}\ol c_{\Si r_{N-1}}) +
\ol\dl(\op\sum_{0\leq |\Si|, |\Xi|}H^{(A,\Xi)(r_{N-1},\Si)}\ol
s_{\Xi A}\ol c_{\Si r_{N-1}})=0,  \label{v145}
\eeq
which is an $(N+1)$-stage Noether identity. Conversely, let
\be
\Phi= \op\sum_{0\leq|\La|} G^{r_N,\La}\ol c_{\La r_N}\om \in
\ol\cP^{0,n}_\infty\{N\}_{N+2}
\ee
be a graded density such that the condition (\ref{v145}) holds.
Then this condition can be extended to a cycle one as follows. It
is brought into the form
\be
&& \dl_N(\op\sum_{0\leq|\La|}  G^{r_N,\La}\ol c_{\La r_N} +
\op\sum_{0\leq |\Si|, |\Xi|}H^{(A,\Xi)(r_{N-1},\Si)}\ol
s_{\Xi A}\ol c_{\Si r_{N-1}})=\\
&& \qquad  -\op\sum_{0\leq|\La|} G^{r_N,\La}d_\La h_{r_N}
+\op\sum_{0\leq |\Si|, |\Xi|}H^{(A,\Xi)(r_{N-1},\Si)}\ol s_{\Xi
A}d_\Si \Delta_{r_{N-1}}.
\ee
A glance at the expression (\ref{v92'}) shows that the term in the
right-hand side of this equality belongs to
$\ol\cP^{0,n}_\infty\{N-2\}_{N+1}$. It is a $\dl_{N-2}$-cycle and,
consequently, a $\dl_{N-1}$-boundary $\dl_{N-1}\Psi$ in accordance
with the $(N+1)$-homology regularity condition. Then the
$(N+1)$-stage Noether identity (\ref{v145}) is a $\ol c_{\Si
r_{N-1}}$-dependent part of the cycle condition
\be
\dl_N(\op\sum_{0\leq|\La|} && G^{r_N,\La}\ol c_{\La r_N} +
\op\sum_{0\leq |\Si|, |\Xi|}H^{(A,\Xi)(r_{N-1},\Si)}\ol s_{\Xi
A}\ol c_{\Si r_{N-1}} -\Psi)=0,
\ee
but $\dl_N\Psi$ does not make a contribution to this identity.

(ii) Being a cycle condition, the $(N+1)$-stage Noether
identity(\ref{v145}) is trivial either if a cycle $\Phi$
(\ref{v156}) is a $\dl_N$-boundary or its summand $G$ is
$\ol\dl$-exact. The $(N+2)$-homology regularity condition implies
that any $\dl_{N-1}$-cycle $\Phi\in
\ol\cP_\infty^{0,n}\{N-1\}_{N+2}\subset
\ol\cP_\infty^{0,n}\{N\}_{N+2}$ is a $\dl_N$-boundary. Therefore,
if $\Phi$ (\ref{v156}) is a representative of a non-trivial
element of $H_{N+2}(\dl_N)$, its summand $G$ linear in $\ol c_{\La
r_N}$ does not vanish. Moreover, it is not a $\ol\dl$-boundary.
Indeed, if $G=\ol\dl \Psi$, then
\mar{v172}\beq
\Phi=\dl_N\Psi +(\ol \dl-\dl_N)\Psi + H. \label{v172}
\eeq
The cycle condition takes the form
\be
\dl_N\Phi=\dl_{N-1}((\ol\dl-\dl_N)\Psi + H)=0.
\ee
Hence, $(\ol \dl-\dl_N)\Psi + H$ is $\dl_N$-exact since any
$\dl_{N-1}$-cycle $\f\in \ol\cP_\infty^{0,n}\{N-1\}_{N+2}$ is a
$\dl_N$-boundary. Consequently, $\Phi$ (\ref{v172}) is a boundary.
If the $(N+2)$-homology regularity condition does not hold,
trivial $(N+1)$-stage Noether identities (\ref{v145}) also come
from non-trivial elements of the homology $H_{N+2}(\dl_N)$.

(iii) Let the $(N+1)$-stage Noether identities be finitely
generated. Namely, there exists a projective Grassmann-graded
$C^\infty(X)$-module $\cC_{(N+1)}$ of finite rank
$\cC_{(N+1)}\subset H_{N+2}(\dl_N)$ such that any element $\Phi\in
H_{N+2}(\dl_N)$ factorizes via elements of $\cC_{(N+1)}$ as
\mar{v160,1}\ben
&& \Phi= \op\sum_{0\leq|\Xi|} \Phi^{r_{N+1},\Xi} d_\Xi
\Delta_{r_{N+1}}\om, \qquad \Phi^{r_{N+1},\Xi}\in
\cS^0_\infty[F;Y], \label{v160}\\
&&\Delta_{r_{N+1}}=G_{r_{N+1}}+ h_{r_{N+1}}=\op\sum_{0\leq|\La|}
\Delta_{r_{N+1}}^{r_N,\La}\ol c_{\La r_N} + h_{r_{N+1}},
\label{v161}
\een
where $\{\Delta_{r_{N+1}}\}$ is local basis for $\cC_{(N+1)}$.
Clearly, this factorization is independent of specification of
this local basis. By virtue of the Serre--Swan theorem,
$\cC_{(N+1)}$ is isomorphic to a module of sections of the product
$\ol V^*_{N+1}\op\times_X \ol E^*_{N+1}$, where $\ol V^*_{N+1}$
and $\ol E^*_{N+1}$ are the density-duals of some vector bundles
$V_{N+1}\to X$ and $E_{N+1}\to X$.
 Let us extend the BGDA $\ol\cP^*_\infty\{N\}$ (\ref{v91})
to the  BGDA $\ol\cP^*_\infty\{N+1\}$ possessing the local basis
\be
\{s^A,\ol s_A, \ol c_r, \ol c_{r_1}, \ldots, \ol c_{r_N}, \ol
c_{r_{N+1}}\}, \quad {\rm Ant}[\ol c_{r_{N+1}}]=N+3, \quad [\ol
c_{r_{N+1}}]=([\Delta_{r_{N+1}}]+1){\rm mod}\,2.
\ee
 It is provided
with the nilpotent graded derivation
\be
\dl_{N+1}=\dl_N + \rdr^{r_{N+1}} \Delta_{r_{N+1}}
\ee
 of antifield
number -1. With this graded derivation, the module
$\ol\cP^{0,n}_\infty\{N+1\}_{\leq N+4}$ of graded densities of
antifield number Ant$[\f]\leq N+4$ is split into the chain complex
\mar{v171}\ben
&&0\lto \im \ol\dl \llr^{\ol\dl} \cP^{0,n}_\infty[\ol Y^*;F;Y;\ol
F^*]_1\llr^{\dl_0} \ol\cP^{0,n}_\infty\{0\}_2\llr^{\dl_1}
\ol\cP^{0,n}_\infty\{1\}_3\cdots
 \label{v171}\\
&& \quad \llr^{\dl_{N-1}} \ol\cP^{0,n}_\infty\{N-1\}_{N+1}
 \llr^{\dl_N} \ol\cP^{0,n}_\infty\{N\}_{N+2}\llr^{\dl_{N+1}}
\ol\cP^{0,n}_\infty\{N+1\}_{N+3}\llr^{\dl_{N+1}}
\ol\cP^{0,n}_\infty\{N+1\}_{N+4}. \nonumber
\een
It is readily observed that this complex is $(N+2)$-exact. In this
case, the $(N+1)$-stage Noether identities (\ref{v145}) come from
the complete $(N+1)$-stage Noether identities
\mar{v162}\beq
 \op\sum_{0\leq|\La|} \Delta_{r_{N+1}}^{r_N,\La} d_\La \Delta_{r_N}\om
+\ol\dl h_{r_{N+1}}\om =0, \label{v162}
\eeq
which are reproduced as the nilpotency conditions of the graded
derivation $\dl_{N+1}$.
\end{proof}

It may happen that the iteration procedure based on Theorem
\ref{v163} is infinite. We restrict our consideration to the case
of a finitely ($N$-stage) reducible Lagrangian system possessing
the finite $(N+2)$-exact chain complex, called the Koszul--Tate
complex,
\mar{w3,4}\ben
&&0\lto \im \ol\dl \llr^{\ol\dl} \cP^{0,n}_\infty[\ol Y^*;F;Y;\ol
F^*]_1\llr^{\dl_0} \ol\cP^{0,n}_\infty\{0\}_2\llr^{\dl_1}
\ol\cP^{0,n}_\infty\{1\}_3\cdots
\label{w3}\\
&& \qquad
 \llr^{\dl_{N-1}} \ol\cP^{0,n}_\infty\{N-1\}_{N+1}
\llr^{\dl_N} \ol\cP^{0,n}_\infty\{N\}_{N+2}\llr^{\dl_N}
\ol\cP^{0,n}_\infty\{N\}_{N+3}, \nonumber\\
&&\dl_N=\rdr^A \cE_A + \op\sum_{0\leq|\La|} \rdr^r
\Delta_r^{A,\La}\ol s_{\La A} + \op\sum_{1\leq k\leq N}\rdr^{r_k}
\Delta_{r_k}, \label{w4}
\een
where $\Delta_{r_k}$ (\ref{v92'}) and the corresponding equalities
(\ref{v93}) are the complete $k$-stage Noether identities. The

\bigskip
\bigskip
\centerline{\bf 3. Noether's inverse second theorem}
\bigskip

Given the BGDA $\ol\cP^*_\infty\{N\}$ (\ref{v91}), let us consider
the BGDA
\mar{w5}\beq
\cP^*_\infty\{N\}=\cP^*_\infty[V_N\cdots V_1V;F;Y;EE_1\cdots E_N]
\label{w5}
\eeq
possessing the local basis
\be
\{s^A, c^r, c^{r_1}, \ldots, c^{r_N}\}, \qquad [c^{r_k}]=([\ol
c_{r_k}]+1){\rm mod}\,2, \qquad  {\rm Ant}[c^{r_k}]=-(k+1),
\ee
and the BGDA
\mar{w6}\beq
P^*_\infty\{N\}=\cP^*_\infty[\ol E^*_N\cdots\ol E^*_1\ol E^*\ol
Y^*V_N\cdots V_1V;F;Y;EE_1\cdots E_N\ol F^*\ol V^*\ol
V^*_1\cdots\ol V_N^*] \label{w6}
\eeq
with the local basis
\mar{w7}\beq
\{s^A, c^r, c^{r_1}, \ldots, c^{r_N},\ol s_A,\ol c_r, \ol c_{r_1},
\ldots, \ol c_{r_N}\}.\label{w7}
\eeq
Following the physical terminology \cite{barn}, we agree to call
$c^{r_k}$, $k\in\Bbb N$, the ghosts of ghost number
gh$[c^{r_k}]=k+1$. Clearly, the BGDAs $\ol\cP^*_\infty\{N\}$
(\ref{v91}) and $\cP^*_\infty\{N\}$ (\ref{w5}) are subalgebras of
the BGDA $P^*_\infty\{N\}$ (\ref{w6}). The Koszul--Tate
differential $\dl_N$ (\ref{w4}) is naturally extended to a graded
derivation of the BGDA $P^*_\infty\{N\}$ (\ref{w6}).

\begin{theo} \label{w35} \mar{w35}
With the Koszul--Tate complex (\ref{w3}), the graded commutative
ring $\cP_\infty^0\{N\}\subset \cP_\infty^*\{N\}$ (\ref{w5}) is
split into the cochain sequence
\mar{w36,108}\ben
&& 0\to \cS^0_\infty[F;Y]\ar^{u_e} \cP^0_\infty\{N\}_1\ar^{u_e}
\cP^0_\infty\{N\}_2\ar^{u_e}\cdots, \label{w36} \\
&& u_e=u + \op\sum_{1\leq k\leq N} u_{(k)},
\label{w108}
\een
graded in a ghost number, where $u$ (\ref{w33}), $u_{(1)}$
(\ref{w38'}) and $u_{(k)}$ (\ref{w38}), $k=2,\ldots, N$,  are the
 gauge, first-stage and higher-stage gauge supersymmetries of an original
Grassmann-graded Lagrangian.
\end{theo}

\begin{proof}
Let us extend an original graded Lagrangian $L$ to
the even graded density
\mar{w8}\beq
L_e=\cL_e\om=L+L_1=L + \op\sum_{0\leq k\leq N} c^{r_k}\Delta_{r_k}\om=L
+\dl_N( \op\sum_{0\leq k\leq N} c^{r_k}\ol c_{r_k}\om), \label{w8}
\eeq
whose summand $L_1$ is linear in ghosts and possesses the zero
antifield number. It is readily observed that $\dl_N(L_e)=0$,
i.e., $\dl_N$ is a variational supersymmetry of the graded
Lagrangian $L_e$ (\ref{w8}). Using the formulas (\ref{0606b}) --
(\ref{0606c}), we obtain
\mar{w16}\ben
 && [\frac{\op\dl^\lto \cL_e}{\dl \ol
s_A}\cE_A +\op\sum_{0\leq k\leq N} \frac{\op\dl^\lto \cL_e}{\dl \ol
c_{r_k}}\Delta_{r_k}]\om = [\frac{\op\dl^\lto \cL_e}{\dl \ol
s_A}\cE_A +\op\sum_{0\leq k\leq N} \frac{\op\dl^\lto \cL_e}{\dl \ol
c_{r_k}} \frac{\dl \cL_e}{\dl
c^{r_k}}]\om= \label{w16}\\
&& \qquad [\up^A\cE_A + \op\sum_{0\leq k\leq N}\up^{r_k}\frac{\dl
\cL_e}{\dl c^{r_k}}]\om= d_H\si, \nonumber\\
&& \up^A= \frac{\op\dl^\lto \cL_e}{\dl \ol s_A}=u^A+w^A
=\op\sum_{0\leq|\La|} c^r_\La\eta(\Delta^A_r)^\La + \op\sum_{1\leq
i\leq N}\op\sum_{0\leq|\La|}
c^{r_i}_\La\eta(\op\dr^\lto{}^A(h_{r_i}))^\La, \nonumber\\
&& \up^{r_k}=\frac{\op\dl^\lto \cL_e}{\dl \ol c_{r_k}} =u^{r_k}+
w^{r_k}= \op\sum_{0\leq|\La|}
c^{r_{k+1}}_\La\eta(\Delta^{r_k}_{r_{k+1}})^\La
+\op\sum_{k+1<i\leq N} \op\sum_{0\leq|\La|}
c^{r_i}_\La\eta(\op\dr^\lto{}^{r_k}(h_{r_i}))^\La. \nonumber
\een
The equality (\ref{w16}) falls into the set of equalities
\mar{w19,',20}\ben
&& \frac{\op\dl^\lto (c^r\Delta_r)}{\dl \ol s_A}\cE_A\om
=u^A\cE_A\om=d_H\si_0, \label{w19}\\
&&  [\frac{\op\dl^\lto (c^{r_1}\Delta_{r_1})}{\dl \ol s_A}\cE_A
+\frac{\op\dl^\lto (c^{r_1}\Delta_{r_1})}{\dl \ol
c_r}\Delta_r]\om= d_H\si_1, \label{w19'}\\
&&  [\frac{\op\dl^\lto (c^{r_i}\Delta_{r_i})}{\dl \ol s_A}\cE_A
+\op\sum_{k<i} \frac{\op\dl^\lto (c^{r_i}\Delta_{r_i})}{\dl \ol
c_{r_k}}\Delta_{r_k}]\om= d_H\si_i, \qquad i=2,\ldots,N.
\label{w20}
\een

A glance at the equality (\ref{w19}) shows that, by virtue of the
decomposition (\ref{g107}), the graded derivation
\mar{w33}\beq
u= u^A\frac{\dr}{\dr s^A}, \qquad u^A
=\op\sum_{0\leq|\La|} c^r_\La\eta(\Delta^A_r)^\La, \label{w33}
\eeq
is a variational supersymmetry of an original graded Lagrangian
$L$. Parameterized by ghosts $c^r$, it is a gauge supersymmetry of
$L$ \cite{jmp05,cmp04}.

The equality (\ref{w19'}) takes the form
\be
&&[\frac{\op\dl^\lto}{\dl \ol s_A}(c^{r_1}h_{r_1}^{(B,\Si)(A,\Xi)}
\ol s_{\Si B}\ol s_{\Xi A})\cE_A + \frac{\op\dl^\lto}{\dl \ol
c_r}(c^{r_1}\op\sum_{0\leq|\Si|}\Delta_{r_1}^{r,\Si}\ol c_{\Si
r})\op\sum_{0\leq|\Xi|} \Delta_r^{B,\Xi}\ol s_{\Xi B}]\om= \\
&& \qquad [\op\sum_{0\leq|\Xi|}
(-1)^{|\Xi|}d_\Xi(c^{r_1}\op\sum_{0\leq|\Si|} 2h_{r_1}^{(B,\Si)(A,\Xi)}
\ol s_{\Si B})\cE_A + u^r\op\sum_{0\leq|\Xi|}
\Delta_r^{B,\Xi}\ol s_{\Xi B}]\om= d_H\si'_1.
\ee
Using the relation (\ref{0606a}), we obtain
\be
[\op\sum_{0\leq|\Xi|} c^{r_1}\op\sum_{0\leq|\Si|} 2h_{r_1}^{(B,\Si)(A,\Xi)} \ol
s_{\Si B} d_\Xi\cE_A + u^r\op\sum_{0\leq|\Xi|} \Delta_r^{B,\Xi}\ol s_{\Xi
B}]\om= d_H\si_1.
\ee
The variational derivative of the both sides of this
equality with respect to the antifield $\ol s_B$ leads to the relation
\be
\op\sum_{0\leq|\Si|} \eta(h_{r_1}^{(B)(A,\Xi)})^\Si d_\Si(2c^{r_1}
d_\Xi\cE_A) + \op\sum_{0\leq|\Si|} u^r_\Si\eta (\Delta^B_r)^\Si=0,
\ee
which is brought into the form
\mar{w34'}\beq
\op\sum_{0\leq|\Si|} d_\Si u^r\frac{\dr}{\dr c^r_\Si} u^B=\ol\dl(\al^B),
\qquad \al^B = -\op\sum_{0\leq|\Si|}
\eta(2h_{r_1}^{(B)(A,\Xi)})^\Si d_\Si(c^{r_1} \ol s_{\Xi A}).
\label{w34'}
\eeq
Therefore, the odd graded derivation
\mar{w38'}\beq
u_{(1)}= u^r\frac{\dr}{\dr c^r}, \qquad u^r=\op\sum_{0\leq|\La|}
c^{r_1}_\La\eta(\Delta^r_{r_1})^\La, \label{w38'}
\eeq
is the first-stage gauge supersymmetry of a
reducible Lagrangian  system \cite{jmp05}.

Every equality (\ref{w20}) is split into a set of equalities with
respect to the polynomial degree in antifields. Let us consider
the one, linear in antifields $\ol c_{r_{i-2}}$ and their jets. We
have
\be
&& [\frac{\op\dl^\lto}{\dl \ol
s_A}(c^{r_i}\op\sum_{0\leq|\Si|,|\Xi|}h_{r_i}^{(r_{i-2},\Si)(A,\Xi)} \ol
c_{\Si r_{i-2}}\ol s_{\Xi A})\cE_A + \\
&& \qquad \frac{\op\dl^\lto}{\dl \ol
c_{r_{i-1}}}(c^{r_i}\op\sum_{0\leq|\Si|}\Delta_{r_i}^{r'_{i-1},\Si}\ol
c_{\Si r'_{i-1}})\op\sum_{0\leq|\Xi|} \Delta_{r_{i-1}}^{r_{i-2},\Xi}\ol
c_{\Xi r_{i-2}}]\om= d_H\si_i.
\ee
It is brought into the form
\be
[\op\sum_{0\leq|\Xi|} (-1)^{|\Xi|}d_\Xi(c^{r_i}\op\sum_{0\leq|\Si|}
h_{r_i}^{(r_{i-2},\Si)(A,\Xi)} \ol c_{\Si r_{i-2}})\cE_A +
u^{r_{i-1}}\op\sum_{0\leq|\Xi|} \Delta_{r_{i-1}}^{r_{i-2},\Xi}\ol c_{\Xi
r_{i-2}}]\om= d_H\si_i.
\ee
Using the relation (\ref{0606a}), we obtain
\be
[\op\sum_{0\leq|\Xi|} c^{r_i}\op\sum_{0\leq|\Si|} h_{r_i}^{(r_{i-2},\Si)(A,\Xi)} \ol
c_{\Si r_{i-2}} d_\Xi\cE_A + u^{r_{i-1}}\op\sum_{0\leq|\Xi|}
\Delta_{r_{i-1}}^{r_{i-2},\Xi}\ol c_{\Xi r_{i-2}}]\om= d_H\si'_i.
\ee
The variational derivative of the both sides of this equality
with respect to the antifield $\ol c_{r_{i-2}}$ leads to the relation
\be
\op\sum_{0\leq|\Si|} \eta(h_{r_i}^{(r_{i-2})(A,\Xi)})^\Si d_\Si(c^{r_i}
d_\Xi\cE_A) + \op\sum_{0\leq|\Si|} u^{r_{i-1}}_\Si\eta
(\Delta^{r_{i-2}}_{r_{i-1}})^\Si=0,
\ee
which takes the form
\mar{w34}\beq
\op\sum_{0\leq|\Si|} d_\Si u^{r_{i-1}}\frac{\dr}{\dr c^{r_{i-1}}_\Si} u^{r_{i-2}}
=\ol\dl(\al^{r_{i-2}}),
\qquad \al^{r_{i-2}} = -\op\sum_{0\leq|\Si|} \eta(h_{r_i}^{(r_{i-2})(A,\Xi)})^\Si
d_\Si(c^{r_i} \ol s_{\Xi A}). \label{w34}
\eeq
Therefore, the odd graded derivations
\mar{w38}\beq
u_{(k)}= u^{r_{k-1}}\frac{\dr}{\dr c^{r_{k-1}}}, \qquad
u^{r_{k-1}}=\op\sum_{0\leq|\La|}
c^{r_k}_\La\eta(\Delta^{r_{k-1}}_{r_k})^\La, \qquad k=2,\ldots,N,
\label{w38}
\eeq
are the $k$-stage gauge supersymmetries \cite{jmp05}. The graded
derivations $u$ (\ref{w33}), $u_{(1)}$ (\ref{w38'}), $u_{(k)}$
(\ref{w38})  are assembled into the ascent operator (\ref{w108})
of ghost number 1. It provides the cochain sequence (\ref{w36}).
\end{proof}

\begin{rem} \label{gg6} \mar{gg6}
The ascent operator (\ref{w108}) need not be nilpotent. We say
that gauge and higher-stage gauge supersymmetries of a Lagrangian
system form an algebra on the shell if the graded derivation
(\ref{w108}) can be extended to a graded derivation $\up$ of ghost
number 1 by means of terms of higher polynomial degree in ghosts
such that $\up$ is nilpotent on the shell. Namely, we have
\mar{w109}\beq
\up=u_e+ \xi= u^A\dr_A + \op\sum_{1\leq k\leq N}(u^{r_{k-1}}
+\xi^{r_{k-1}})\dr_{r_{k-1}}+\xi^{r_N}\dr_{r_N}, \label{w109}
\eeq
where all the coefficients $\xi^{r_{k-1}}$ are at least quadratic
in ghosts and $(\up\circ\up)(f)$ is $\ol\dl$-exact for any graded
function $f\in \cP^0_\infty\{N\}\subset P^0_\infty\{N\}$. This
nilpotency condition falls into a set of equalities with respect
to the polynomial degree in ghosts. Let us write the first and
second of them involving the coefficients $\xi_2^{r_{k-1}}$
quadratic in ghosts. We have
\mar{w110,3}\ben
&& \op\sum_{0\leq|\Si|} d_\Si u^r\dr^\Si_r u^B=\ol\dl(\al^B_1), \qquad
\op\sum_{0\leq|\Si|} d_\Si u^{r_{k-1}}\dr_{r_{k-1}}^\Si u^{r_{k-2}}
=\ol\dl(\al^{r_{k-2}}_1), \quad 2\leq k\leq N, \label{w110}\\
&&  \op\sum_{0\leq|\Si|}[d_\Si u^A\dr^\Si_A u^B
+d_\Si\xi^r_2\dr^\Si_r u^B]=\ol\dl(\al^B_2),
\label{w111} \\
&& \op\sum_{0\leq|\Si|}[d_\Si u^A\dr^\Si_A u^{r_{k-1}} +
d_\Si\xi^{r_k}_2\dr^\Si_{r_k} u^{r_{k-1}} +
d_\Si
u^{r'_{k-1}}\dr^\Si_{r'_{k-1}}\xi^{r_{k-1}}_2]= \ol\dl(\al^{r_{k-1}}_2),
 \label{w112}\\
&& \xi_2^r=\xi^{r,\La,\Si}_{r',r''} c^{r'}_\La c^{r''}_\Si,
\qquad  \xi_2^{r_k}=\xi^{r_k,\La,\Si}_{r,r'_k} c^r_\La c^{r'_k}_\Si, \qquad
2\leq k\leq N. \label{w113}
\een
The equalities (\ref{w110}) reproduce the relations (\ref{w34'}) and
(\ref{w34}) in Theorem \ref{w35}. The equalities (\ref{w111}) -- (\ref{w112})
provide the generalized
commutation relations on the shell between gauge and higher-stage gauge
supersymmetries, and one  can think of the coefficients $\xi_2$ (\ref{w113})
as being {\it sui generis} generalized structure functions
\cite{jmp05,fulp02}.
\end{rem}

\bigskip
\bigskip

\centerline{\bf 4. Example}

\bigskip

We address the topological BF theory of two exterior forms $A$ and
$B$ of form degree $|A|+|B|=\di X-1$ on a smooth manifold $X$
\cite{birm}, but restrict our consideration to its simplest
variant where $A$ is a function \cite{jpa05,jmp05a}.

Let us consider the fiber bundle
\be
Y=\Bbb R\op\times_X \op\w^{n-1} T^*X,
\ee
coordinated by $(x^\la, A, B_{\m_1\ldots \m_{n-1}})$ and provided with
the canonical $(n-1)$-form
\be
B=\frac{1}{(n-1)!}B_{\m_1\ldots \m_{n-1}}dx^{\m_1}\w\cdots\w
dx^{\m_{n-1}}.
\ee
The Lagrangian and the Euler--Lagrange operator of the topological
BF theory  read
\mar{v182,3}\ben
&& L_{\rm BF}=\frac1n Ad_HB, \label{v182}\\
&&\dl L= dA\w \cE\om + dB_{\m_1\ldots \m_{n-1}}\w \cE^{\m_1\ldots
\m_{n-1}}\om, \nonumber\\
&& \cE=\e^{\m\m_1\ldots \m_{n-1}} d_\m B_{\m_1\ldots \m_{n-1}},
\qquad \cE^{\m_1\ldots \m_{n-1}} = - \e^{\m\m_1\ldots
\m_{n-1}}d_\m A, \label{v183}
\een
where $\e$ is the Levi--Civita symbol. Let consider the BGDA
$\cP^*_\infty[\ol Y^*;Y]$ where
\be
VY=Y\op\times_X Y, \qquad \ol Y^*= (\Bbb R\op\times_X
\op\w^{n-1}TX)\op\ot_X \op\w^n T^*X.
\ee
It  possesses the local basis $\{ A, B_{\m_1\ldots \m_{n-1}}, \ol
s, \ol s^{\m_1\ldots \m_{n-1}}\}$, where $\ol s, \ol s^{\m_1\ldots
\m_{n-1}}$ are odd antifields of antifield number 1. With the
nilpotent Koszul--Tate differential
\be
\ol\dl=\frac{\rdr}{\dr \ol s}\cE + \frac{\rdr}{\dr \ol
s^{\m_1\ldots \m_{n-1}}} \cE^{\m_1\ldots \m_{n-1}},
\ee
we have the complex (\ref{v042}):
\be
0\lto \im\ol\dl \llr^{\ol\dl} \cP^{0,n}_\infty[\ol Y^*;Y]_1
\llr^{\ol\dl} \cP^{0,n}_\infty[\ol Y^*;Y]_2.
\ee
A generic one-chain reads
\be
\Phi= \op\sum_{0\leq |\La|}(\Phi^\La\ol s_\La +
\Phi^\La_{\m_1\ldots \m_{n-1}} \ol s^{\m_1\ldots
\m_{n-1}}_\La)\om,
\ee
and the cycle condition takes the form
\mar{v189}\beq
\ol\dl\Phi=\Phi^\La\cE_\La + \Phi^\La_{\m_1\ldots \m_{n-1}}
\cE^{\m_1\ldots \m_{n-1}}_\La=0. \label{v189}
\eeq
If $\Phi^\La$ and $\Phi^\La_{\m_1\ldots\m_{n-1}}$ are independent
of the variational derivatives (\ref{v183}) (i.e., $\Phi$ is a
nontrivial cycle), the equality (\ref{v189}) is split into the
following ones
\be
\Phi^\La\cE_\La=0, \qquad
 \Phi^\La_{\m_1\ldots \m_{n-1}} \cE^{\m_1\ldots
\m_{n-1}}_\La=0.
\ee
The first equality holds iff $\Phi^\La=0$, i.e., there is
no Noether identity involving $\cE$. The second one is
satisfied iff
\be
\Phi^{\la_1\ldots \la_k}_{\m_1\ldots\m_{n-1}}\e^{\m\m_1\ldots
\m_{n-1}}=- \Phi^{\m\la_2\ldots
\la_k}_{\m_1\ldots\m_{n-1}}\e^{\la_1\m_1\ldots \m_{n-1}}.
\ee
It follows that $\Phi$ factorizes as
\be
\Phi= \op\sum_{0\leq |\Xi|} G_{\nu_2\ldots\nu_{n-1}}^\Xi
d_\Xi\Delta^{\nu_2\ldots\nu_{n-1}}\om
\ee
via local graded densities
\mar{v190}\beq
\Delta^{\nu_2\ldots\nu_{n-1}}=\Delta^{\nu_2\ldots\nu_{n-1},
\la}_{\al_1\ldots\al_{n-1}}\ol
s^{\al_1\ldots\al_{n-1}}_\la=\dl^\la_{\al_1}\dl^{\nu_2}_{\al_2}\cdots
\dl^{\nu_{n-1}}_{\al_{n-1}} \ol s^{\al_1\ldots\al_{n-1}}_\la=
d_{\nu_1}\ol s^{\nu_1\nu_2\ldots\nu_{n-1}}, \label{v190}
\eeq
which provide the complete Noether identities
\mar{v191}\beq
d_{\nu_1}\cE^{\nu_1\nu_2\ldots\nu_{n-1}}=0. \label{v191}
\eeq

The local graded densities (\ref{v190}) form the basis for a
projective $C^\infty(X)$-module of finite rank which is isomorphic
to the module of sections of the vector bundle
\be
\ol V^*=\op\w^{n-2} TX\op\ot_X \op\w^n T^*X, \qquad V= \op\w^{n-2}
T^*X.
\ee
Therefore, let us enlarge the BGDA $\cP^*_\infty[\ol Y^*;Y]$ to
the BGDA $\ol\cP^*_\infty\{0\}= \cP^*_\infty[\ol Y^*;Y;V]$
possessing the local basis $\{A, B_{\m_1\ldots \m_{n-1}}, \ol s,
\ol s^{\m_1\ldots \m_{n-1}}, \ol c^{\m_2\ldots \m_{n-1}}\}$, where
$\ol c^{\m_2\ldots \m_{n-1}}$ are even antifields of antifield
number 2. We have the nilpotent graded derivation
\be
\dl_0= \ol\dl + \frac{\rdr}{\dr \ol c^{\m_2\ldots \m_{n-1}}}
\Delta^{\m_2\ldots \m_{n-1}}
\ee
of $\cP^*_\infty\{0\}$. Its nilpotency is equivalent to the
complete Noether identities (\ref{v191}). Then we obtain the
one-exact complex
\be
0\lto \im\ol\dl \llr^{\ol\dl} \cP^{0,n}_\infty[\ol Y^*;Y]_1
\llr^{\dl_0} \ol\cP^{0,n}_\infty\{0\}_2 \llr^{\dl_0}
\ol\cP^{0,n}_\infty\{0\}_3.
\ee

Iterating the arguments, we come to the $(N+1)$-exact complex
(\ref{v94}) for $N\leq n-3$ as follows. Let us consider the
corresponding BGDA
\be
\ol\cP^*_\infty\{N\}= \cP^*_\infty[...V_3V_1\ol Y^*;Y;VV_2V_4...],
\qquad V_k=\op\w^{n-k-2} T^*X, \qquad k=1,\ldots, N,
\ee
possessing the local basis
\be
&&\{A, B_{\m_1\ldots \m_{n-1}}, \ol s, \ol s^{\m_1\ldots
\m_{n-1}}, \ol c^{\m_2\ldots \m_{n-1}},\ldots,\ol c^{\m_{N+2}\ldots \m_{n-1}}\},\\
&& [\ol c^{\m_{k+2}\ldots \m_{n-1}}]=(k+1){\rm mod}\,2, \qquad
{\rm Ant}[\ol c^{\m_{k+2}\ldots \m_{n-1}}]=k+3.
\ee
It is provided with the nilpotent graded derivation
\mar{va202}\beq
\dl_N=\dl_0 + \op\sum_{1\leq k\leq N}\frac{\rdr}{\dr \ol
c^{\m_{k+2}\ldots \m_{n-1}}}
 \Delta^{\m_{k+2}\ldots \m_{n-1}}, \qquad
\Delta^{\m_{k+2}\ldots
\m_{n-1}}=d_{\m_{k+1}}\ol c^{\m_{k+1}\m_{k+2}\ldots \m_{n-1}},
\label{va202}
\eeq
of antifield number -1. Its nilpotency results from the Noether
identities (\ref{v191}) and equalities
\mar{v212}\beq
d_{\m_{k+2}}\Delta^{\m_{k+2}\ldots \m_{n-1}}=0, \qquad k\in \Bbb
N, \label{v212}
\eeq
which are $k$-stage Noether identities \cite{jpa05}. Then the
manifested $(N+1)$-exact complex reads
\mar{v203}\ben
&&0\lto \im \ol\dl \llr^{\ol\dl} \cP^{0,n}_\infty[\ol
Y^*;Y]_1\llr^{\dl_0} \ol\cP^{0,n}_\infty\{0\}_2\llr^{\dl_1}
\ol\cP^{0,n}_\infty\{1\}_3\cdots
\label{v203}\\
&& \qquad
 \llr^{\dl_{N-1}} \ol\cP^{0,n}_\infty\{N-1\}_{N+1}
\llr^{\dl_N} \ol\cP^{0,n}_\infty\{N\}_{N+2}\llr^{\dl_N}
\ol\cP^{0,n}_\infty\{N\}_{N+3}. \nonumber
\een
It obeys the following $(N+2)$-homology regularity condition.

\begin{lem} \label{v220} \mar{v220}
Any $(N+2)$-cycle $\Phi\in \ol\cP^{0,n}_\infty\{N-1\}_{N+2}$ up to
a $\dl_{N-1}$-boundary is
\mar{v218}\ben
&& \Phi=\op\sum_{k_1+\cdots +k_i+3i=N+2}\sum_{0\leq|\La_1|,\ldots,
|\La_i|}G^{\La_1\cdots \La_i}_{\m^1_{k_1+2}\ldots
\m^1_{n-1};\ldots; \m^i_{k_i+2}\ldots \m^i_{n-1}} \label{v218} \\
&&\qquad d_{\La_1} \Delta^{\m^1_{k_1+2}\ldots \m^1_{n-1}}\cdots
d_{\La_i} \Delta^{\m^i_{k_i+2}\ldots \m^i_{n-1}}\om, \qquad
k=-1,0,1,\ldots, N, \nonumber
\een
where $\ol c^{\m_1\ldots\m_{n-1}}=\ol s^{\m_1\ldots\m_{n-1}}$ and
$\Delta^{\m_1\ldots\m_{n-1}}=\cE^{\m_1\ldots\m_{n-1}}$. It follows
that $\Phi$ is a $\dl_N$-boundary.
\end{lem}

\begin{proof}
Let us choose some basis element $\ol c^{\m_{k+2}\ldots \m_{n-1}}$
and denote it simply by $\ol c$. Let $\Phi$ contain a summand
$\f_1 \ol c$, linear in $\ol c$. Then the cycle condition reads
\be
\dl_{N-1}\Phi=\dl_{N-1}(\Phi-\f_1 \ol c) + (-1)^{[\ol
c]}\dl_{N-1}(\f_1)\ol c + \f \Delta=0, \qquad \Delta=\dl_{N-1}\ol
c.
\ee
It follows that $\Phi$ contains a summand $\psi\Delta$ such that
\be
(-1)^{[\ol c]+1}\dl_{N-1}(\psi)\Delta +\f\Delta=0.
\ee
This equality implies the relation
\mar{v213}\beq
\f_1=(-1)^{[\ol c]+1}\dl_{N-1}(\psi) \label{v213}
\eeq
because the reduction conditions (\ref{v212}) involve total
derivatives of $\Delta$, but not $\Delta$. Hence,
\be
\Phi=\Phi' +\dl_{N-1}(\psi \ol c),
\ee
where $\Phi'$ contains no term linear in $\ol c$. Furthermore, let
$\ol c$ be even and $\Phi$ has a summand $\sum \f_r \ol c^r$
polynomial in $\ol c$. Then the cycle condition leads to the
equalities
\be
\f_r\Delta=-\dl_{N-1}\f_{r-1}, \qquad r\geq 2.
\ee
Since $\f_1$ (\ref{v213}) is $\dl_{N-1}$-exact, then $\f_2=0$ and,
consequently, $\f_{r>2}=0$. Thus, a cycle $\Phi$ up to a
$\dl_{N-1}$-boundary contains no term polynomial in $\ol c$. It
reads
\mar{v217}\beq
\Phi=\op\sum_{k_1+\cdots +k_i+3i=N+2}\sum_{0<|\La_1|,\ldots,
|\La_i|}G^{\La_1\cdots \La_i}_{\m^1_{k_1+2}\ldots
\m^1_{n-1};\ldots; \m^i_{k_i+2}\ldots \m^i_{n-1}} \ol
c^{\m^1_{k_1+2}\ldots \m^1_{n-1}}_{\La_1}\cdots \ol
c_{\La_i}^{\m^i_{k_i+2}\ldots \m^i_{n-1}}\om. \label{v217}
\eeq
However, the terms polynomial in $\ol c$ may appear under general
covariant transformations
\be
\ol c'^{\nu_{k+2}\ldots \nu_{n-1}}=\det(\frac{\dr x^\al}{\dr
x'^\bt}) \frac{\dr x'^{\nu_{k+2}}}{\dr x^{\m_{k+2}}}\cdots
\frac{\dr x'^{\nu_{n-1}}}{\dr x^{\m_{n-1}}}\ol c^{\m_{k+2}\ldots
\m_{n-1}}
\ee
of a chain $\Phi$ (\ref{v217}). In particular, $\Phi$ contains the
summand
\be
\op\sum_{k_1+\cdots +k_i+3i=N+2}F_{\nu^1_{k_1+2}\ldots
\nu^1_{n-1};\ldots; \nu^i_{k_i+2}\ldots \nu^i_{n-1}} \ol
c'^{\nu^1_{k_1+2}\ldots \nu^1_{n-1}}\cdots \ol
c'^{\nu^i_{k_i+2}\ldots \nu^i_{n-1}},
\ee
which must vanish if $\Phi$ is a cycle. This takes place only if
$\Phi$ factorizes through the graded densities
$\Delta^{\m_{k+2}\ldots \m_{n-1}}$ (\ref{va202}) in accordance
with the expression (\ref{v218}).
\end{proof}

Following the proof of Lemma \ref{v220}, one can also show that
any $(N+2)$-cycle $\Phi\in \ol\cP^{0,n}_\infty\{N\}_{N+2}$ up to a
boundary takes the form
\be
\Phi=\op\sum_{0\leq|\La|}G^\La_{\m_{N+2}\ldots \m_{n-1}} d_\La
\Delta^{\m_{N+2}\ldots \m_{n-1}}\om,
\ee
i.e., the homology $H_{N+2}(\dl_N)$ of the complex (\ref{v203}) is
finitely generated by the cycles $\Delta^{\m_{N+2}\ldots
\m_{n-1}}$. Thus, the complex (\ref{v203}) admits the
$(N+2)$-exact extension (\ref{v171}).

The iteration procedure is prolonged till $N=n-3$. We have the
BGDA $\ol\cP^*\{n-2\}$, where $V_{n-2}=X\times\Bbb R$. It
possesses the local basis
\be
\{A, B_{\m_1\ldots \m_{n-1}}, \ol s, \ol s^{\m_1\ldots \m_{n-1}},
\ol c^{\m_2\ldots \m_{n-1}},\ldots,\ol c^{\m_{n-1}}, \ol c\},
\ee
where $[\ol c]=(n-1){\rm mod}\,2$ and Ant$[\ol c]=n+1$.
The corresponding Koszul--Tate complex reads
\be
&&0\lto \im \ol\dl \llr^{\ol\dl} \cP^{0,n}_\infty[\ol
Y^*;Y]_1\llr^{\dl_0} \ol\cP^{0,n}_\infty\{0\}_2\llr^{\dl_1}
\ol\cP^{0,n}_\infty\{1\}_3\cdots\\
&& \qquad
 \llr^{\dl_{n-3}} \ol\cP^{0,n}_\infty\{n-3\}_{n-1}
\llr^{\dl_{n-2}} \ol\cP^{0,n}_\infty\{n-2\}_n\llr^{\dl_{n-2}}
\ol\cP^{0,n}_\infty\{n-2\}_{n+1}. \\
&&\dl_{n-2}=\dl_0 + \op\sum_{1\leq k\leq n-3}\frac{\rdr}{\dr \ol
c^{\m_{k+2}\ldots \m_{n-1}}} \Delta^{\m_{k+2}\ldots \m_{n-1}} +
\frac{\rdr}{\dr \ol c}\Delta, \qquad \Delta=d_{\m_{n-1}}\ol
c^{\m_{n-1}}.
\ee

Let us enlarge the BGDA $\ol\cP^*_\infty\{n-2\}$ to the BGDA
$P^*_\infty\{n-2\}$ (\ref{w6}) with the local basis
\be
\{A, B_{\m_1\ldots \m_{n-1}}, c_{\m_2\ldots \m_{n-1}},\ldots,
c_{\m_{n-1}},  c,
 \ol s, \ol s^{\m_1\ldots \m_{n-1}},
\ol c^{\m_2\ldots \m_{n-1}},\ldots,\ol c^{\m_{n-1}}, \ol c\},
\ee
where $c_{\m_2\ldots \m_{n-1}},\ldots, c_{\m_{n-1}},  c$ are the
corresponding ghosts, and let us consider the BGDA
$\cP^*_\infty\{n-2\}$ with the local basis $\{A, B_{\m_1\ldots
\m_{n-1}}, c_{\m_2\ldots \m_{n-1}},\ldots, c_{\m_{n-1}}, c\}$. By
virtue of Theorem \ref{w35}, the graded commutative ring
$\cP^0_\infty\{n-2\}$ is split into the cochain sequence
\mar{gg10,4}\ben
&& 0\to \cO^0_\infty Y\ar^{u_e} \cP^0_\infty\{n-2\}_1\ar^{u_e}
\cP^0_\infty\{n-2\}_2\ar^{u_e}\cdots, \label{gg10} \\
&& u_e=u + \op\sum_{1\leq k\leq n-2} u_{(k)}, \label{gg11}\\
&& u=-d_{\m_1}c_{\m_2\ldots\m_{n-1}}\frac{\dr}{\dr
B_{\m_1\ldots\m_{n-1}}}, \label{gg12}\\
&& u_{(k)}=-d_{\m_{k+1}}c_{\m_{k+2}\ldots\m_{n-1}}\frac{\dr}{\dr
c_{\m_{k+1}\ldots\m_{n-1}}}, \qquad 1\leq k\leq n-3,
\label{gg13}\\
&& u_{(n-2)}= -d_\m c \frac{\dr}{\dr c_\m}, \label{gg14}
\een
where $u$ (\ref{gg12}) and $u_{(k)}$ (\ref{gg13}) -- (\ref{gg14})
are the gauge and higher-stage gauge supersymmetries of the
Lagrangian (\ref{v182}) \cite{jpa05}. It is readily observed that
the ascent operator (\ref{gg11}) is nilpotent, i.e., the sequence
(\ref{gg10}) is a cochain complex.

\bigskip
\bigskip
\centerline{\bf 5. Appendix}
\bigskip

The proof of Theorem \ref{v11} follows that of \cite{cmp04},
Theorem 2.1 when $Y$ is an affine bundle.

\begin{lem} \label{v38'} \mar{v38'}
If $Y=\Bbb R^{n+m}\to \Bbb R^n$, the complex (\ref{g111}) at all
the terms, except $\Bbb R$, is exact, while the complex
(\ref{g112}) is exact.
\end{lem}

\begin{proof}
This is the case of an affine bundle $Y$, and the above mentioned
exactness has been proved when the ring $\cO^0_\infty Y$ is
restricted to the subring $\cP^0_\infty Y$ of polynomial functions
(see \cite{cmp04}, Lemmas 4.2 -- 4.3). The proof of these lemmas
is straightforwardly extended to $\cO^0_\infty Y$ if the homotopy
operator (4.5) in \cite{cmp04}, Lemma 4.2 is replaced with that
(4.8) in \cite{cmp04}, Remark 4.1.
\end{proof}

We first prove Theorem \ref{v11} for the above mentioned BGDA
$\G(\gQ^*_\infty[F;Y])$. Similarly to $\cS^*_\infty[F;Y]$, the
sheaf $\gQ^*_\infty[F;Y]$ and the BGDA $\G(\gQ^*_\infty[F;Y])$ are
split into the variational bicomplexes, and we consider their
subcomplexes
\mar{v35-8}\ben
&& 0\ar \Bbb R\ar \gQ^0_\infty[F;Y]\ar^{d_H}\gQ^{0,1}_\infty[F;Y]
\cdots \ar^{d_H} \gQ^{0,n}_\infty[F;Y]\ar^\dl {\got E}_1,
\label{v35}\\
&& 0\to \gQ^{1,0}_\infty[F;Y]\ar^{d_H} \gQ^{1,1}_\infty[F;Y]\cdots
\ar^{d_H}\gQ^{1,n}_\infty[F;Y]\ar^\vr {\got E}_1\to 0, \label{v36}\\
&& 0\ar \Bbb R\ar
\G(\gQ^0_\infty[F;Y])\ar^{d_H}\G(\gQ^{0,1}_\infty[F;Y]) \cdots
\ar^{d_H} \G(\gQ^{0,n}_\infty[F;Y])\ar^\dl \G({\got E}_1),
\label{v37} \\
&&  0\to \G(\gQ^{1,0}_\infty[F;Y])\ar^{d_H}
\G(\gQ^{1,1}_\infty[F;Y])\cdots
\ar^{d_H}\G(\gQ^{1,n}_\infty[F;Y])\ar^\vr \G({\got E}_1)\to 0,
\label{v38}
\een
where ${\got E}_1 =\vr(\gQ^{1,n}_\infty[F;Y])$. By virtue of Lemma
\ref{v38'}, the complexes (\ref{v35}) -- (\ref{v36}) at all the
terms, except $\Bbb R$, are exact. The terms
$\gQ^{*,*}_\infty[F;Y]$ of the complexes (\ref{v35}) --
(\ref{v36}) are sheaves of $\G(\gQ^0_\infty)$-modules. Since
$J^\infty Y$ admits a partition of unity just by elements of
$\G(\gQ^0_\infty)$, these sheaves are fine and, consequently,
acyclic. By virtue of the abstract de Rham theorem (see
\cite{cmp04}, Theorem 8.4, generalizing \cite{hir}, Theorem
2.12.1), cohomology of the complex (\ref{v37}) equals the
cohomology of $J^\infty Y$ with coefficients in the constant sheaf
$\Bbb R$ and, consequently, the de Rham cohomology of $Y$, which
is the strong deformation retract of $J^\infty Y$. Similarly, the
complex (\ref{v38}) is proved to be exact. It remains to prove
that cohomology of the complexes (\ref{g111}) -- (\ref{g112})
equals that of the complexes (\ref{v37}) -- (\ref{v38}). The proof
of this fact straightforwardly follows the proof of \cite{cmp04},
Theorem 2.1, and it is a slight modification of the proof of
\cite{cmp04}, Theorem 4.1, where graded exterior forms on the
infinite order jet manifold $J^\infty Y$ of an affine bundle are
treated as those on $X$.

{\sc \hfill Depart. of Math. and Inform.

\hfill University of Camerino

\hfill 62032 Camerino (MC), Italy

\medskip

\hfill Depart. of Theor. Phys.

\hfill Moscow State University

\hfill 117234 Moscow, Russia }

\end{document}